%% file: H_PrincQuad.tex
\documentclass[reqno]{amsart}
\usepackage{etex}
\usepackage{amsfonts}
\usepackage{amscd}
\usepackage{amsbsy}
\usepackage{mathtools}
\usepackage{amsxtra}
\usepackage{amssymb}
\usepackage{epsfig}
\usepackage{epic,eepic}
\usepackage{graphicx}
\usepackage[all]{xy}
\usepackage{xypic}
\usepackage{psfrag}
\usepackage{amsmath}
\usepackage{amsthm}
\usepackage{setspace}
\usepackage{hyperref}
\usepackage{url}
\usepackage{here}
\usepackage{multirow}
\usepackage{todonotes}
\newlength{\halfbls}
\usepackage{tikz}
\usepackage{tikz}
\usetikzlibrary{calc}
\usetikzlibrary{decorations.pathreplacing,decorations.markings}
\usetikzlibrary{arrows}
\usetikzlibrary{decorations.markings, decorations.shapes}
\usetikzlibrary{arrows}
\usetikzlibrary{positioning}
\usetikzlibrary{intersections}
\usetikzlibrary{decorations.pathmorphing}

\usepackage{thmtools}

\DeclareRobustCommand{\SkipTocEntry}[9]{}

\usepackage{filecontents}
\include{macros_QuadRec}

\makeatletter
\renewcommand{\@bibtitlestyle}{%
  \@xp\part\@xp*\@xp{\refname}%
}
\makeatother

\makeatletter
\newcommand{\addresseshere}{%
  \enddoc@text\let\enddoc@text\relax
}
\makeatother

\makeatletter \let\@wraptoccontribs\wraptoccontribs
\def\@setauthors{%
  \begingroup
  \def\thanks{\protect\thanks@warning}%
  \trivlist
  \centering\footnotesize \@topsep30\p@\relax
  \advance\@topsep by -\baselineskip
  \item\relax
  \author@andify\authors
  \def\\{\protect\linebreak}%
  \MakeUppercase{\authors}%
  \ifx\@empty\contribs
  \else
    \linebreak \@setcontribs
    \@closetoccontribs
  \fi
  \endtrivlist
  \endgroup
}
\makeatother

\title[Volumes of strata of quadratic differentials]
{Masur-Veech volumes and intersection theory:
the principal strata of quadratic differentials}

\author{Dawei Chen}
\thanks{Research of D.C. is partially supported by the NSF CAREER grant DMS-1350396 and a von Neumann Fellowship at the IAS in Spring 2019.}
\address{Department of Mathematics, Boston College, Chestnut Hill, MA 02467, USA}
\email{dawei.chen@bc.edu}

\author{Martin M\"oller}
\thanks{Research of M.M. is partially supported  
  by the DFG-project MO 1884/1-1, by the LOEWE-Schwerpunkt
  ``Uniformisierte Strukturen in Arithmetik und Geometrie'',
  and by the NSF grant DMS-1440140 while the author was in residence at
  the MSRI, Berkeley during Fall 2019.
}
\address{
Institut f\"ur Mathematik, Goethe--Universit\"at Frankfurt,
Robert-Mayer-Str. 6--8,
60325 Frankfurt am Main, Germany
}
\email{moeller@math.uni-frankfurt.de}

\author{Adrien Sauvaget}
\thanks{Research of A.S. is supported by the Dutch Research Council (NWO) grant 613.001.651.}
\address{
Mathematical Institute, Utrecht University, Budapestlaan 6 /
Hans Freudenthal Bldg, 3584 CD Utrecht, The Netherlands
}
\email{a.c.b.sauvaget@uu.nl}

\contrib[with an appendix by]{Ga\"etan Borot}
\address{Max-Planck-Institut f\"ur Mathematik,
Vivatsgasse 7, 
53111 Bonn, Germany
}
\email{gborot@mpim-bonn.mpg.de}
\thanks{The work of G.B., A.G., and D.L.\ benefits from the support of
  the Max-Planck-Gesellschaft}
\contrib[]{Alessandro Giacchetto}
\address{
  Max-Planck-Institut f\"ur Mathematik,
  Vivatsgasse 7, 
53111 Bonn, Germany
}
\email{agiacche@mpim-bonn.mpg.de}

\contrib[]{Danilo Lewa\'nski}
\address{
  Max-Planck-Institut f\"ur Mathematik,
  Vivatsgasse 7, 
53111 Bonn, Germany
}
\email{danilo.lewanski@studio.unibo.it}
      
\begin{document}

\maketitle
\begin{abstract}
We describe a conjectural formula via intersection numbers for the Masur-Veech volumes of strata of quadratic differentials with prescribed zero orders, and we prove the formula for the case when the zero orders are odd. For the principal strata of quadratic differentials with simple zeros, the formula reduces to compute the top Segre class of the quadratic Hodge bundle, which can be further simplified to certain linear Hodge integrals.  An appendix proves that the intersection of this class with $\psi$-classes can be computed by Eynard-Orantin topological recursion.
\par
As applications, we analyze numerical properties of Masur-Veech volumes, area Siegel-Veech constants and sums of Lyapunov exponents of the principal strata for fixed genus and varying number of zeros, which settles the corresponding conjectures due to Grivaux-Hubert, Fougeron, and elaborated in~\cite{TRMVV}. We also describe conjectural formulas for area Siegel-Veech constants and sums of Lyapunov exponents for arbitrary affine invariant submanifolds, and verify them for the principal strata.  
\end{abstract}
\tableofcontents
\SaveTocDepth{1} 

\input{sec_intro} 

\input{sec_IntMV}


\input{sec_principal}


\input{sec_svLyap}


\SaveTocDepth{0} 
\def\appendixname{}
\appendix
\def\contribs{Ga\"etan Borot, Alessandro Giacchetto,
and Danilo Lewa\'nski}
\def\appendixtitle{A second topological recursion for Masur--Veech volumes}
\markboth{\MakeUppercase\contribs}{\MakeUppercase\appendixtitle}

%
%
%
%
%
%

\makeatletter
\def\part{\@startsection{part}{0}%
\z@{\linespacing\@plus\linespacing}{.5\linespacing}%
  {\normalfont\Large\bfseries\raggedright}}
\makeatother
\part*{A second topological recursion for Masur--Veech volumes}
\noindent by \textsc{Ga\"etan Borot, Alessandro Giacchetto,}
and \textsc{Danilo Lewa\'nski}


\vspace{0.5cm}
\renewcommand \thesection {\Alph{section}}
\input{sec_appendix}


\renewcommand{\bibliofont}{}

\newcommand{\etalchar}[1]{$^{#1}$}
\def\cprime{$'$}

\end{document}

%% file: macros_QuadRec.tex
\DeclareMathAlphabet\mathbfcal{OMS}{cmsy}{b}{n}
\usepackage{xparse}

\theoremstyle{plain}
\newtheorem{Defi}{Definition}[section]  
   \newtheorem{Prop}[Defi]{Proposition}
\newtheorem{Lemma}[Defi]{Lemma}    \newtheorem{Cor}[Defi]{Corollary}
\newtheorem{Thm}[Defi]{Theorem}

\newtheorem{conjecture}[Defi]{Conjecture}

\theoremstyle{remark}

\newtheorem{Rem}[Defi]{Remark}

\newcommand{\CC}{\mathbb{C}}

\newcommand{\BM}{\overline{\mathcal M}}
\newcommand{\BQ}{\overline{\mathcal Q}}

\newcommand{\BNN}{\overline{\mathcal{N}}}
\newcommand{\NN}{\mathbb{N}}
 
\newcommand{\PP}{\mathbb{P}} 
\newcommand{\QQ}{\mathbb{Q}} 
\newcommand{\RR}{\mathbb{R}}

\newcommand{\ZZ}{\mathbb{Z}}

 \newcommand{\cEE}{{\mathcal E}} 
 \newcommand{\cXX}{{\mathcal X}}  
\newcommand{\cCC}{{\mathcal C}} 
 
\newcommand{\cHH}{{\mathcal H}} 
\newcommand{\cMM}{{\mathcal M}} \newcommand{\cOO}{{\mathcal O}}
 \newcommand{\cNN}{{\mathcal N}}
 
\newcommand{\cQQ}{{\mathcal Q}}

\newcommand{\area}{{\rm area}} 

\newcommand{\irr}{\operatorname{irr}}
\DeclareMathOperator{\vol}{vol}

\newcommand{\ch}{\operatorname{ch}}

\newcommand{\Area}{\operatorname{Area}}
\newcommand{\wt}{\widetilde}

\def\={\;=\;}  \def\+{\,+\,}

\DeclareMathAlphabet{\eucal}{U}{eus}{m}{n}
\DeclareMathAlphabet{\newcal}{U}{dutchcal}{m}{n}

   \def\wt#1{\widetilde{#1}}   \def\wh#1{\widehat{#1}}

\def\GL{{\textrm{GL}}}

\def\be{\begin{equation}}   \def\ee{\end{equation}}     \def\bes{\begin{equation*}}    \def\ees{\end{equation*}}
\def\ba{\be\begin{aligned}} \def\ea{\end{aligned}\ee}   \def\bas{\bes\begin{aligned}}  \def\eas{\end{aligned}\ees}

\newcommand{\moduli}[1][g]{{\mathcal M}_{#1}}
\newcommand{\omoduli}[1][g]{{\Omega \mathcal M}_{#1}}

\newcommand{\barmoduli}[1][g]{{\overline{\mathcal M}}_{#1}}

\newcommand{\barcQ}[1][g]{{\overline{\mathcal Q}}_{#1}}

\newcommand{\oM}{\overline{\mathcal M}}

\newcommand{\BE}{\overline{\mathbb E}}

\newcounter{savedtocdepth}
\newcommand*{\SaveTocDepth}[1]{%
  \addtocontents{toc}{%
    \protect\setcounter{savedtocdepth}{\protect\value{tocdepth}}%
    \protect\setcounter{tocdepth}{#1}%
  }%
}

\DeclareDocumentCommand{\LMS}{ O{\mu} O{g}} {\Xi\mathcal{M}_{#2}(#1)}

%% file: sec_intro.tex
\section{Introduction}
\label{sec:intro}

Computing volumes of moduli spaces via intersection theory has significance in many aspects.  For example, the Weil-Petersson volumes of the moduli spaces of marked Riemann surfaces can be calculated by intersection numbers on the Deligne-Mumford compactification $\BM_{g,n}$, which motivated Mirzakhani to give a proof of Witten's conjecture by hyperbolic geometry \cite{MirzWP}.  A more recent instance is the intersection-theoretic formula for the Masur-Veech volumes of moduli spaces of Abelian differentials with prescribed orders of zeros, which can be used to determine the large genus asymptotics of the Masur-Veech volumes \cite{CMSZ}.  
\par
The moduli spaces of Abelian differentials carry a natural $\GL_2^+(\mathbb R)$-action induced by varying the translation surface representations of Abelian differentials. 
 The orbit closures of this action are called affine invariant submanifolds, since they have locally linear structures \cite{esmi, esmimo}.  Besides the ambient moduli spaces, affine invariant submanifolds can provide interesting (and more challenging) playgrounds for us to detect relevant geometric invariants, such as volumes and intersection numbers.  A prominent type of affine invariant submanifolds arises from moduli spaces of quadratic differentials (also called half-translation surfaces), which can be lifted into the corresponding moduli spaces of Abelian differentials via the canonical double cover.  In this paper we focus on the moduli spaces of quadratic differentials.  
 \par
\subsection*{Masur-Veech volumes}
Let $(\mu, \nu)$ be an integer partition of $4g-4$, where $\mu = (2m_i)_{i=1}^r$
are the even parts and $\nu = (2n_j-1)_{j=1}^s$ are the odd parts, with $m_i, n_j \geq 0$. Consider the moduli space (also called the stratum) $\cQQ_{g,r+s}(\mu, \nu)$  parameterizing quadratic differentials $q$ on Riemann surfaces of genus $g$ such that $q$ has $r$ even order zeros of type $\mu$ and $s$ odd order zeros of type $\nu$. Note that $q$ is allowed to have simple poles, i.e., when some $n_j = 0$, which are regarded as ``zeros of order $-1$'' because in this case the surface still has finite area under the metric of $q$. Similarly $q$ is allowed to have ordinary marked points, i.e., when some $m_i = 0$, which are regarded as ``zeros of order $0$''.  
Moreover, in our setting all zeros are labeled (as in \cite{CMSZ}, but contrary to \cite{DGZZ} where only the simple poles are labeled).  
\par 
Let $\PP\cQQ_{g,r+s}(\mu, \nu) = \cQQ_{g,r+s}(\mu, \nu)/\mathbb C^{*}$ be the projectivized stratum parameterizing quadratic differentials of type $(\mu, \nu)$ up to scaling.  Denote by $\PP\barcQ[g,r+s](\mu,\nu)$ its closure in the incidence variety
compactification (IVC) of the strata of quadratic differentials \cite{kdiff}.  Let $\zeta$ be
the first Chern class of the universal line bundle $\cOO(1)$ on $\PP\barcQ[g,r+s](\mu,\nu)$. Denote by 
$\psi_{i}$ the cotangent line bundle class on $\BM_{g,n}$ associated with the $i$-th marked point as well as its pullbacks to the strata $\PP\BQ_{g,r+s}(\mu, \nu)$. We first describe a conjectural formula to compute the Masur-Veech volumes of all strata $\cQQ_{g,r+s}(\mu, \nu)$ via intersection theory (see Section~\ref{sec:MVnorm} for our convention on volume normalization). 
\par
\begin{conjecture}\label{conj:main}
The Masur-Veech volumes of strata of quadratic differentials can be obtained as the intersection numbers
\ba
\label{eq:conj-main}
\vol(\cQQ_{g,r+s}(\mu, \nu)) &\=  \frac{2^{r-s+3}(2\pi {\rm i})^{2g-2+s}}{(2g-3+r+s)!}
\int_{\PP\barcQ[g,r+s](\mu,\nu)}\zeta^{2g+s-3} \psi_{1} \cdots \psi_{r}\,,
\ea
where $\psi_1, \ldots, \psi_r$ are associated with the $r$ even order zeros. 
\end{conjecture}
\par 
Note that for the strata of quadratic differentials, the number $r$ of even order zeros is equal to the dimension of
the relative period foliation. Hence from the viewpoint of period coordinates, when there is no relative period, i.e., when $r = 0$, the volume formula of the strata of quadratic differentials behaves similarly to that of the minimal strata $\cHH(2g-2)$ of Abelian differentials \cite{SauvagetMinimal}. We can thus prove the above conjectural formula for this special case using a good metric on the tautological line bundle $\cOO(-1)$ \cite{CoMoZa}. 
\par 
\begin{Thm} 
\label{thm:conjtrue}
For the strata of quadratic differentials with odd zeros only, we have 
 \ba
 \label{eq:odd}
\vol(\cQQ_{g, s}(\nu)) &\=  \frac{2^{3-s}(2\pi {\rm i})^{2g-2+s}}{(2g-3+s)!}
\int_{\PP\barcQ[g,s](\nu)}\zeta^{2g+s-3}\,. 
\ea
\end{Thm}
\par
We remark that we have also verified the conjectural formula for a number of low genus strata that have relative periods, i.e. when $r\neq 0$, by ad hoc calculations. However, to prove the formula in full generality, one either needs a good metric on the $\psi$-bundles or a volume recursion out of merging zeros, which we plan to study in future work.
\par  
Theorem~\ref{thm:conjtrue} is particularly useful for the principal strata of quadratic differentials with only simple zeros (and simple poles). Let $\barcQ[g,n]$ be the quadratic Hodge bundle over $\BM_{g,n}$ whose fiber over a stable pointed curve $(X, p_1, \ldots, p_n)$ is $H^0(\omega_X^{\otimes 2} (p_1 + \cdots + p_n))\cong \mathbb C^{3g-3+n}$, where $\omega_X$ is the dualizing line bundle of $X$. The interior space $\cQQ_{g,n}$ over $\cMM_{g,n}$ parameterizes quadratic differentials with at worst simple poles at the marked points. Hence $\barcQ[g,n]$ provides an alternate compactification (smaller than the IVC in \cite{kdiff}) for the principal stratum $\cQQ_{g,n}(1^{4g-4+n}, -1^{n})$ (here the $4g-4+n$ simple zeros are not labeled as they can merge to form higher order zeros in the quadratic Hodge bundle). In this case the top self-intersection of the $\zeta$-class corresponds to the top Segre class $s(\barcQ[g,n])$ of the quadratic Hodge bundle  $\barcQ[g,n]$. Characteristic classes of the (Abelian) Hodge bundle and its variants including the $k$-th Hodge bundle were computed in \cite{mumford83, Chiodo} by the Grothendieck-Riemann-Roch formula. Combining their results with intersection calculations on $\BM_{g,n}$, 
Formula~\eqref{eq:odd} can then be evaluated more explicitly in terms of linear Hodge integrals as follows.  
\par 
\begin{Thm} 
\label{thm:PSINT}
The Masur-Veech volumes of the principal strata of quadratic differentials can
be obtained as the intersection numbers  
\ba
\label{eq:principal}
\, &\phantom{\=} \,\,\vol(\cQQ_{g,4g-4+2n}(1^{4g-4+n},-1^n)) \\
\, &\= \frac{2^{2g+1}\pi^{6g-6+2n}}{(6g-7+2n)!}  \sum_{i=0}^g \frac{(4g-4+n)!}{(2g-3+n+i)!} 
\int_{\BM_{g, 2g-3+2n+i}} \!\!\!\!\!\!\!\!\psi_{n+1}^2 \cdots \psi_{2g-3+2n+i}^2 \lambda_{g-i}\,  \\ 
\, &\=  \frac{2^{2g+1}\pi^{6g-6+2n}}{(6g-7+2n)!}  \sum_{i=0}^g \frac{(4g-4+n)!(4g-7+2n+i)!!}{(2g-3+i)! (4g-7 + i)!!} 
\int_{\BM_{g, 2g-3+i}} \!\!\!\!\!\!\!\!\psi_{1}^2 \cdots \psi_{2g-3+i}^2 \lambda_{g-i}\,,  
\ea
where $\lambda_i$ is the $i$-th Chern class of the (Abelian) Hodge bundle on $\BM_{g,k}$.   
\end{Thm}
\par
In addition, the Appendix shows that the intersection of this Segre class with $\psi$-classes
\begin{equation}
\label{Segrepsi} \int_{\BM_{g,n}} s(\barcQ[g,n]) \prod_{i = 1}^n \psi_i^{k_i}\qquad (k_i \geq 0)
\end{equation} 
can be computed by the Eynard-Orantin topological recursion for the spectral curve
\begin{equation}
\label{spsp}\left\{\begin{array}{lll} x(z) & = & -z - \ln z \\ y(z) & = & z^2 \end{array}\right.\qquad \omega_{0,2}(z_1,z_2) \= \frac{{\rm d}z_1{\rm d}z_2}{(z_1 - z_2)^2}\,.
\end{equation}
In particular, the volumes of the principal strata can be recovered from the $k_i = 0$ term. In~\cite{TRMVV} another set of numbers $F_{g,n}(k_1,\ldots,k_n)$ was constructed and computed by the topological recursion, which can also be expressed via intersection theory on $\barmoduli[g,n]$ such that the $k_i = 0$ case recovers the volumes of the principal strata. Rather surprisingly though, in \cite{TRMVV} the spectral curve is very different from \eqref{spsp} and these $F_{g,n}$ are not a priori related to \eqref{Segrepsi} except for $k_i = 0$, where it gives a different expression of $\vol(\cQQ_{g,4g - 4 + 2n}(1^{4g - 4 + n},-1^{n}))$ as the top intersection of a class on $\barmoduli[g,n]$.
\par
Theorem~\ref{thm:PSINT} can be used to analyze numerical properties of the Masur-Veech volumes of the principal strata with fixed genus and varying number of zeros.  
\par 
 \begin{Thm}[{\cite[Conjecture 5.4 and (5.12)]{TRMVV}}] 
 \label{thm:the7conj54}
For all $g\geq 1$ and for $n\geq 0$ (except for $g =1$ and $n \leq 1$ where the strata are empty), there exist two rational polynomials $p_g(n)$
and $q_g(n)$ of degree $\lfloor (g-1)/2 \rfloor$ and $\lfloor g/2 \rfloor$
respectively, such that 
\ba
\label{eq:volume-fix-g}
\, &\phantom{\=} \,\,\frac{\vol(\cQQ_{g,4g-4+2n}(1^{4g-4+n}, -1^n))}{\pi^{6g-6+2n}} \\
\, & \= 2^n  \frac{(2g-3+n)! (4g-4+n)!}{(6g-7+2n)!} \big(p_g(n)+ \gamma_{2g-3+n} q_g(n)\big)\,,
 \ea
where $ \gamma_k= \frac{1}{4^k}\binom{2k}{k}$. 
\par 
Moreover for fixed $g$ and $n\to \infty$, we have the asymptotic growth rates 
 \ba
 \label{eq:volume-large-n}
\vol(\cQQ_{g,4g-4+2n}(1^{4g-4+n}, -1^n))\, \sim \, 2^{-n} \pi^{6g-6+2n + \epsilon (g)/2} m_g n^{g/2}\,, 
\ea
where $\epsilon (g) = 0$ or $1$ is the parity of $g$ and where $2^{6g-7} m_g $
is the top coefficient of $q_g$ if $g$ is even or the top coefficient of $p_g$
if $g$ is odd. 
 \end{Thm}
 \par
 For the case of tori, the above results can be
 described more explicitly as follows.  
 \par 
\begin{Cor}
\label{cor:volume-g=1}
 For $g=1$ and $n\geq 2$, we have 
 \ba
 \label{eq:volume-g=1}
 \vol(\cQQ_{1,2n}(1^n, -1^n)) & \=  \pi^{2n} \frac{n!}{3 (2n-1)!} \big((2n-3)!! + (2n-2)!!\big)\,. \\
 \ea
 \end{Cor}
\par
We remark that the coefficient $m_g$ in Theorem~\ref{thm:the7conj54} is a rescaling of the intersection number $\int_{\oM_{g,3g-3}} \psi_1^2 \ldots \psi_{3g-3}^2$ and can be computed efficiently (see the end of Section~\ref{subsec:volume-large-n} for references on this topic).  In addition to the large~$n$ asymptotic, in~\cite{DGZZ} the large $g$ asymptotic of $\vol(\cQQ_{g,n})$ was described conjecturally, and the conjecture was further extended and refined in~\cite{YZZVolume}. 
\par 
\subsection*{Lyapunov exponents and Siegel-Veech constants} We now switch gears to discuss applications of our results in surface dynamics.  
Consider the straight line flow on a torus with a
half-translation structure induced by a quadratic differential with~$\ell$ simple zeros and~$\ell$ simple poles.
Closing up a random trajectory as it comes within the $(1/n)$-ball
of its starting point defines a collection of cycles~$\{\gamma_n\}_{n \in \NN}$.
The logarithm of the size of $\gamma_n$ in any norm on the cohomology
of the torus tends to $\boldsymbol{\lambda}_1$ times the logarithm of the flat
length of~$\gamma_n$, where $\boldsymbol{\lambda}_1$ is a quantity that does
not depend on the starting point and the direction of the trajectory,
as long as they are generic. This is a consequence of Oseledets theorem
and~$\boldsymbol{\lambda}_1$ is the (first) Lyapunov exponent of the straight line flow
(see e.g.~\cite{zoSQ, zorich06} for more details about Lyapunov exponents). 
\par
Near each of the simple poles, the trajectory makes a $U$-turn. Pulling the
trajectory tight as a cohomology class should cause 
drastic shortcuts, and hence the growth rate of the cohomology classes
$\gamma_n$ is expected to decrease with the number of poles~$\ell$. This conjecture
was first announced by Grivaux and Hubert. In the case of half-translation surfaces
in the stratum~$\cQQ(\ell, -1^\ell)$ with a zero of order~$\ell$
and $\ell$~simple poles, Fougeron~\cite{FouLNP} found an upper bound for
$\boldsymbol{\lambda}_1$ in the order
of~$1/\ell$. His method works in any sequence of strata where the
largest order of zeros is unbounded. In the case of simple zeros,
Fougeron conjectured a decay in the order of~$1/\sqrt{\ell}$. A refined version of 
these conjectures for the principal strata of quadratic differentials appeared in~\cite{TRMVV} in terms of area Siegel-Veech constants $c_{\area}$, which can determine the sums of (involution-invariant) Lyapunov exponents $L^{+}$, and vice versa by \cite[Theorem 2]{ekz}. These conjectures were stated as conditional corollaries in~\cite{TRMVV} assuming the numerical results in our Theorem~\ref{thm:the7conj54}. Here we prove them unconditionally based on the following formulas that express $c_{\area}$ and $L^+$ as intersection numbers.   
\par
\begin{Thm}
\label{thm:sv-L}
For the principal strata of quadratic differentials, we have 
\ba
\label{eq:sv-principal}
c_{\area}(1^{4g-4+n}, -1^n) & \= -\frac{1}{2\pi^2} \frac{ \int_{\PP\barcQ[g,n]}\zeta^{6g-8+2n}\delta }{ \int_{\PP\barcQ[g,n]}\zeta^{6g-7+2n}}\,,
\ea
where $\delta$ is the total boundary divisor class, and   
\ba
\label{eq:L-principal}
L^+(1^{4g-4+n}, -1^n) & \= - 2 \frac{ \int_{\PP\barcQ[g,n]}\zeta^{6g-8+2n}\lambda_1 }{ \int_{\PP\barcQ[g,n]}\zeta^{6g-7+2n}}\,. 
\ea  
\end{Thm}
\par 
 \begin{Cor}[{\cite[Corollary 5.5 and (5.12)]{TRMVV}}] 
 \label{cor:sv-L}
For all $g\geq 1$ and $n\geq 0$ (except for $g =1$ and $n \leq 1$ where the strata are empty), there exist two rational polynomials $p^{*}_g(n)$ and $q^{*}_g(n)$ of degree $\lfloor (g+3)/2 \rfloor$ and $1+\lfloor g/2 \rfloor$ respectively, such that  
\ba
\label{eq:sv-fixed-g}
 c_{\area}(1^{4g-4+n}, -1^n) \= \frac{1}{\pi^2} \frac{\frac{p_g^{*}(n)}{2g-3+n} + \gamma_{2g-3+n}q_g^{*}(n)}{p_g(n) + \gamma_{2g-3+n}q_g(n)}\,,
 \ea
 where $p_g$ and $q_g$ are the polynomials in Theorem~\ref{thm:the7conj54}. 
 \par 
 Moreover, there exist two rational polynomials $r_g(n)$ and $s_g(n)$ of degree $\lfloor g/2 \rfloor$ and $\lfloor (g+1)/2 \rfloor$  respectively, such that 
 \ba
\label{eq:L-fixed-g-large-n}
 L^+(1^{4g-4+n}, -1^{n}) & \=  \frac{1}{2g-3+n} \frac{r_g(n) + \gamma_{2g-3+n}s_g(n)}{p_g(n) + \gamma_{2g-3+n}q_g(n)} \\
 & \,\sim\,  \frac{\pi^{1/2-\epsilon (g)} n_g/m_g}{\sqrt{n}}
\ea
as $n\to \infty$, where $m_g$ is defined in Theorem~\ref{thm:the7conj54} and  
$2^{6g-7} n_g$ is the top coefficient of $r_g$ if $g$ is even or 
the top coefficient of $s_g$ if $g$ is odd.  
 \end{Cor}
\par 
Again for the case of tori, the above results can be described more explicitly. 
\par 
 \begin{Cor}
 \label{cor:sv-L-g=1}
 For $g=1$ and $n\geq 2$, we have 
 \ba
 \label{eq:sv-g=1}
 c_{\area}(1^n, -1^n) & \= \frac{1}{\pi^2} \Bigg(     \frac{n}{6} + \frac{6}{1 + \frac{(2n-2)!!}{(2n-3)!!}} \Bigg)\,
 \ea
 and 
 \ba
  \label{eq:L-g=1}
 L^+(1^n, -1^n) & \=  \frac{2}{1+ \frac{(2n-2)!!}{(2n-3)!!}}\,.
 \ea
 \end{Cor}
 \par
 Motivated by Conjecture~\ref{conj:main} and Theorem~\ref{thm:sv-L}, we come up with an analogous conjecture for $c_{\area}$ and $L^+$ for all strata of quadratic differentials. 
 \par 
 \begin{conjecture}
   \label{conj:sv-L}
   The Siegel-Veech constants and Lyapunov exponents
for the strata $\cQQ_g(\mu, \nu)$ with $\mu = (2m_i)_{i=1}^r$ and $\nu = (2n_j-1)_{j=1}^s$ are given by 
$$ c_{\area}(\mu, \nu) \= -\frac{1}{2\pi^2} \frac{ \int_{\PP\barcQ(\mu,\nu)}\zeta^{2g+s-4} \psi_{1} \cdots \psi_{r} \delta }{ \int_{\PP\barcQ(\mu,\nu)}\zeta^{2g+s-3} \psi_{1} \cdots \psi_{r}}\,, $$
$$ L^+(\mu, \nu) \= - 2\frac{ \int_{\PP\barcQ(\mu,\nu)}\zeta^{2g+s-4} \psi_{1} \cdots \psi_{r} \lambda_1}{ \int_{\PP\barcQ(\mu,\nu)}\zeta^{2g+s-3} \psi_{1} \cdots \psi_{r}}\,, $$
where $\psi_1, \ldots, \psi_r$ are associated with the $r$ even order zeros. 
\end{conjecture}
 \par
We remark that Kontsevich speculated an implicit intersection formula to compute sums of Lyapunov exponents for strata of Abelian differentials~\cite[Section~7]{kontsevich}, and such a formula was justified explicitly in our previous work~\cite{CMSZ}. 
\par
Motivated by the above results, at the end of the paper we make a general conjecture to compute area Siegel-Veech constants and sums of Lyapunov exponents as intersection numbers for an arbitrary affine invariant submanifold (see Conjecture~\ref{conj:ais-sv-L}), and we summarize known cases as evidences of the conjecture.   
 \par
 \subsection*{Related works}  We briefly review several related works about the Masur-Veech volumes of strata of quadratic differentials. A standard method to compute volumes is to determine the quasimodular forms that arise from pillowcase covers and compute their large degree asymptotics (which follows from 
 the initial idea of counting torus covers for Abelian differentials in \cite{eo}). This was used in \cite{GouExplicit} to obtain explicit values of volumes for a number of low dimensional strata of quadratic differentials. Another approach is to decompose half-translation surfaces
into ribbon graphs and sum up the corresponding local contributions. This was first carried out for all strata in genus zero in \cite{aez} and then extended 
to the principal strata for all genus in \cite{DGZZ}.  These local contributions of ribbon graphs are indeed intersection numbers of $\psi$-classes on $\BM_{g,n}$ that come from Kontsevich's proof of Witten's conjecture \cite{KontWitten}. Moreover, in \cite{TRMVV} the same sum of local contributions was shown to arise as constant terms of a family of polynomials that are determined by topological recursion, whose approach relies on statistics of hyperbolic curves. Our method is different from all of these works, as we use a good metric on the compactified strata and simplify the intersection calculation for the principal strata by the Grothendieck-Riemann-Roch formula and linear Hodge integrals.  
\par 
\subsection*{Organization of the paper} In Section~\ref{sec:IntMV} we show that when there is no relative period, the Hermitian metric induced by the area form gives the Masur-Veech volume form up to an explicit scaling factor, thus proving Theorem~\ref{thm:conjtrue}.  In Section~\ref{sec:principal} we reduce the volume formula for the principal strata to certain linear Hodge integrals and analyze its numerical properties for fixed genus and varying number of zeros, proving Theorem~\ref{thm:PSINT}, Theorem~\ref{thm:the7conj54} and Corollary~\ref{cor:volume-g=1}. In Section~\ref{sec:svL} we apply our results to area Siegel-Veech constants and sums of Lyapunov exponents of the principal strata, proving Theorem~\ref{thm:sv-L}, Corollary~\ref{cor:sv-L} and Corollary~\ref{cor:sv-L-g=1}.
The Appendix establishes a topological recursion for \eqref{Segrepsi} and shows how to compute the volumes of the principal strata from it.
\par 
\subsection*{Acknowledgements} We thank Vincent Delecroix, Maksym Fedorchuk, Elise Goujard, Shuai Guo, Felix Janda, Brian Lehmann, Motohico Mulase, Di Yang, Don Zagier, Youjin Zhang, Zhengyu Zong, Anton Zorich, and Dimitri Zvonkine for helpful discussions on related topics.  
Part of the work was carried out when the first two authors attended the CMO workshop ``Flat Surfaces and Dynamics on Moduli Space, II'' in May 2019. We thank the organizers for the invitation and hospitality.  

%% file: sec_IntMV.tex
\section{Masur-Veech volumes as intersection numbers}
\label{sec:IntMV}

In this section we prove the expression of Masur-Veech
volumes as intersection numbers in Theorem~\ref{thm:conjtrue} for the strata of quadratic differentials with no even order zeros. The argument
is largely parallel to the proof of \cite[Proposition~1.3]{SauvagetMinimal}.
This section also serves the purpose of introducing period coordinates
and explaining the volume normalization convention we use.
\par
We first set up some general notation. We write $N = r + s$ for the total number
of marked zeros. Half-translation surfaces parameterized in $\cQQ_{g,N}(\mu,\nu)$ are usually denoted by $(X,q)$. For a surface $(X, q) \in \cQQ_{g,N}(\mu,\nu)$, let $\pi\colon \wh{X} \to X$ be the canonical double cover such that $\pi^*q = \omega^2$ is the square of an Abelian differential (see e.g.\ \cite[Section 2.2]{ekz}). We simply denote by $\wh{\cHH}$ the space of lifts of $\cQQ_{g,N}(\mu,\nu)$ via the double cover, with all preimages of the singularities of $q$ being labeled. This locus is an affine invariant submanifold of the stratum
$\omoduli[\wh g,\wh N](\wh \mu)$ of Abelian differentials, where
$\wh g = 2g -1 + s/2$, $\wh N = 2r+s$, and $\wh \mu = ((m_i,m_i)_{i=1}^r, (n_j)_{j=1}^s)$. Note that the squaring map $\wh{\cHH} \to \cQQ_{g,N}(\mu,\nu)$
is finite of degree $2^{r+1}$, which is due to labeling the $r$ pairs of preimages of the even order zeros as well as choosing the sign of $\omega$. 
Therefore, the induced map  $\PP\wh{\cHH} \to \PP\cQQ_{g,N}(\mu,\nu)$ on
the projectivized spaces has degree $2^r$, since $\pm\omega$ correspond to the same point in the projectivization. 
In particular, $\PP\wh{\cHH} \to \PP\cQQ_{g,N}(\mu,\nu)$ is an isomorphism for the case $r = 0$ (though $\PP\wh{\cHH}$ carries an order-$2$ stacky structure due to the involution of the double cover). 
\par
\subsection{The Masur-Veech volume form} \label{sec:MVnorm}
The Masur-Veech volume form is defined on $\cQQ_{g,N}(\mu,\nu)$ using period
coordinates and gives a finite measure on the hyperboloid $\cQQ_{g,N}^{1}(\mu,\nu)$
of half-translation surfaces of area $1/2$, i.e. the double cover surface $( \wh{X}, \omega)$ has area one. This volume form is obtained by desintegration of the Lebesgue measure with respect to the area coordinate. The volume of $\cQQ_{g,N}^{1}(\mu,\nu)$ can be computed by integration of a volume form on  $\PP\cQQ_{g,N}(\mu,\nu)$ as we describe below.  
\par
We start with the flat area form on the double cover
\bes
h(\omega) \coloneqq \Area_{\wh{X}}(\omega) \= \frac{{\rm i}}{2}
\int_{\wh{X}} \omega \wedge \bar{\omega}
= \frac{{\rm i}}{2} \sum_{i=1}^{\wh g} (z_{A_i}\bar{z}_{B_i} - z_{B_i}\bar{z}_{A_i})\,, 
\ees
where $(A_i, B_i)_{i=1}^{\wh g}$ form a symplectic basis of $H_1(\wh{X},\ZZ)$
and $(z_{A_i}, z_{B_i})_{i=1}^{\wh g}$ are the $\omega$-periods of this basis.
Let $\tau$ be the involution on~$\wh{X}$ whose quotient map is~$\pi$.
\par
Both $\wh{\cHH}$ and $\cQQ_{g,N}(\mu, \nu)$ are locally modeled on the $\tau$-anti-invariant cohomology $H_-^{1}(\wh{X}, \CC)$. Any choice of scale of the Lebesgue measure that is invariant under the symplectic group results in an (infinite) measure $\wt{\nu}_{MV}$ on $\cQQ_{g,N}(\mu,\nu)$. For our choice of normalization we define
\be
H_\pm^{1}(\wh{X}, \ZZ) \= H^{1}(\wh{X}, \ZZ) \cap
H_\pm^{1}(\wh{X}, \RR)\,.
\ee
Consequently, $H_{-}^{1}(\wh{X}, \ZZ \oplus i\ZZ)$ is a lattice in 
$H_-^{1}(\wh{X}, \CC)$ and we normalize $\wt{\nu}_{MV}$ so that this lattice
has covolume one.
\par
With the help of the area normalization and of $\wt{\nu}_{MV}$ we define 
\be
\nu_{MV}(U) \= \wt{\nu}_{MV}(C_U)\,, \quad (U \subseteq \cQQ_{g,N}^{1}(\mu,\nu))\,,
\ee
where $C_U = \{ \lambda (X,q) \colon (X,q) \in U, \lambda \in [0,1]\}$
is the cone under~$U$. This measure is finite by \cite{masur82, veech82}.
Abusing notation we also denote by $\nu_{MV}$ the pushforward of this measure
on $\PP\cQQ_{g,N}(\mu,\nu)$. By the squaring map we can view $\nu_{MV}$
also as a measure on $\PP\wh{\cHH}$. We write
\bes
\vol(\cQQ_{g,N}(\mu,\nu)) \= {\rm dim}_{\RR}(\widehat{\mathcal{H}}) \cdot \nu_{MV}(\PP\cQQ_{g,N}(\mu,\nu))
\ees
for the Masur-Veech volume of the total space.

\subsection{The metric on $\cOO_{\PP\wh{\cHH}}(-1)$ and comparison of volume forms}

The area form~$h$ induces a Hermitian metric (still denoted by $h$) on $\cOO_{\PP\wh{\cHH}}(-1)$. 
For a section $\sigma$ of $\cOO_{\PP\wh{\cHH}}(-1)$, we 
consider the associated curvature $(1,1)$-form 
\bes
\omega_h \= \frac{1}{2\pi {\rm i}} \partial\bar{\partial} \log h (\sigma)
\quad \text{and} \quad
\nu_h \= \omega_h^{2g-3+s}\,
\ees
the corresponding volume form on $\PP\wh{\cHH}$. 
\par
\begin{Lemma} \label{le:VFcomp}
For $r=0$, i.e. when there is no zero of even order, the two volume forms are proportional as follows: 
$$ \nu_h \=
- \frac{(2g-3+s)!}{2^{2g+1}(2\pi {\rm i})^{2g-2+s}}\, \dim_\RR(\wh{\cHH})\, \nu_{MV}\,. $$
\end{Lemma}
\par
\begin{proof}
Fix a point $(\wh X, \omega) \in \wh{\cHH}$. An open neighborhood
of this point can be written under the period coordinates as $\wh X + v$
for~$v \in H_-^1(\wh X, \CC)$ small enough. After shrinking we may assume
that this neighborhood is contained in the positive cone 
$$ \cCC = \{\wh{X} + v \in H^1(\wh X, \CC)\colon  h(\wh X+v) > 0\}\,. $$
It thus suffices to prove the volume form relation on $\PP\cCC$. The
proof follows the idea of~\cite[Lemma~2.1]{SauvagetMinimal}, with an extra twist
coming from the fact that the sum of the two subspaces $H^1_\pm(\wh X, \ZZ)$
spans a proper subgroup of finite index in $H^1(\wh{X}, \ZZ)$. 
\par
First note that since $r = 0$, $s$ is positive and even, and hence the rank of the
anti-invariant part is greater than or equal to~$2g$. By \cite[Corollary~12.1.5]{bl}
the symplectic form restricted to $H_1^+(\wh X, \ZZ)$ is of type $(2,\ldots,2)$ and
the restriction to $H_1^-(\wh X,\ZZ)$ is of type $(2,\ldots,2,1,\ldots1)$
with~$g$ entries of $2$. 
\par
By the Riemann-Hurwitz formula the genus $\wh g$ of $\wh X$ is given by $\wh g = 2g-1+s/2$. We also set 
$\wt{g} = g-1+s/2 = \wh g - g$ to simplify formulas later.  We define  $C_i^\pm = A_{2i-1} \pm A_{2i}$ and
$D_i^\pm = B_{2i-1} \pm B_{2i}$ for $i=1,\ldots,g$. For an appropriate order
of the elements in the symplectic basis we have 
\be
H_1^-(\wh X,\ZZ) = \langle C_1^-,D_1^-,\ldots,C_g^-,D_g^-, A_{2g+1},
B_{2g+1}, \ldots, A_{\wh g}, B_{\wh g} \rangle_{\ZZ}\,,
\ee
see also \cite[Section~2.1]{GouExplicit} for a realization of these cycles. We define $C_{g+i}^- = A_{2g+i}$ and $D_{g+i}^- = B_{2g+i}$ for $1\leq i \leq \wt g - g$, 
define $\chi_j = 1/2$ for $1\leq j \leq g$ and $\chi_j = 1$ for $g< j \leq \wt g$, and denote by $z_{C_j^-}$ and $z_{D_j^-}$ the corresponding $\omega$-periods of $C_j^-$ and $D_j^-$. Since $\omega$ pairs trivially with the invariant eigenspace $H_1^+(\wh X, \ZZ)$, for $1\leq j \leq g$ we have 
$$z_{C_j^-} = 2z_{A_{2j-1}} = -2  z_{A_{2j}}\,, \quad z_{D_j^-} = 2z_{B_{2j-1}} = -2  z_{B_{2j}}\,,$$ 
 and hence 
 $$ \frac{1}{2} (z_{C_j^-}\bar{z}_{D_j^-} - z_{D_j^-}\bar{z}_{C_j^-}) = (z_{A_{2j-1}}\bar{z}_{B_{2j-1}} - z_{B_{2j-1}}\bar{z}_{A_{2j-1}}) + (z_{A_{2j}}\bar{z}_{B_{2j}} - z_{B_{2j}}\bar{z}_{A_{2j}})\,.$$
 It follows that  
\be
h(\omega) \= \frac{{\rm i}}{2} \sum_{j=1}^{\wt{g}}
\chi_j (z_{C_j^-}\bar{z}_{D_j^-} - z_{D_j^-}\bar{z}_{C_j^-})\,.
\ee
Passing to the coordinate system 
\bes z_{c_j} \= \frac{1}{2} (z_{C_j^-} - {\rm i} z_{D_j^-})\,, \quad
z_{d_j} \= \frac{1}{2} (z_{C_j^-} + {\rm i} z_{D_j^-})\, \ees
of $H^1_-(\wh X,\CC)$, the Hermitian metric $h$ (of signature $(\wt{g},
\wt{g})$) can be written as 
\be \label{eq:defmetricgg}
h(\{z_{c_j},z_{d_j}\}_{j=1}^{\wt g}) \= \sum_{j=1}^{\wt{g}}
\chi_j (z_{c_j}\bar{z}_{c_j} - z_{d_j}\bar{z}_{d_j})\,.
\ee
\par
As in the proof of \cite[Lemma~2.1]{SauvagetMinimal}, we now proceed by comparing $\omega_h$ and $\nu_{MV}$ to the forms $\omega'_h$ and $\nu'_{MV}$ obtained from the (positive
definite) metric~$h'$ with a plus sign (instead of minus) in~\eqref{eq:defmetricgg}.
Since the ratios are invariant under the group $U(\wt{g},
\wt{g}) \cap U(2\wt{g})$, it suffices to compare $\nu_h$ and $\nu_{MV}$
on a fundamental domain for this group inside the cone, i.e.\ the
set $(z_{c_1}, z_{d_1},0,\ldots,0)$ in the projectivized cone $\PP \cCC$.
By symmetry consideration it suffices to focus on the chart $U_{c_1} =
\{z_{c_1} = 1\}$ and use the section $\sigma(z_{d_1},z_{c_2},\ldots)
= (1,z_{d_1},z_{c_2},\ldots)$ of $\cOO_{\PP\wh{\cHH}}(-1)$. 
\par
In this chart, we claim that the Masur-Veech volume form is
\be \label{eq:compare}
\nu_{MV} \= \frac{2(2\pi){\rm i}^{2\wt{g}-1}}{\dim_\RR(\wh{\cHH})h(\sigma)^{2\wt{g}}}
\Bigl( dz_{d_1}\wedge d\bar{z}_{d_1} \wedge \prod_{j = 2}^{\wt g}
( dz_{c_j}\wedge d\bar{z}_{c_j} \wedge dz_{d_j}\wedge d\bar{z}_{d_j})
\Bigr)\,.
\ee
To see this, first note that the factor $2\pi$ comes from integrating the argument of 
the (complex) coordinate $z_{c_1}$ with fixed norm. Next, the volume of a cone over a base parameterizing surfaces of area $h$ is $h^{\dim_\CC(\wh{\cHH})}$ times 
the volume of the corresponding cone over the base parameterizing surfaces of area one, thus explaining the factor $h(\sigma)^{2\wt{g}}$. Moreover, one checks that 
$$dz_{c_j} \wedge d\bar{z}_{c_j} \wedge dz_{d_j} \wedge d\bar{z}_{d_j} = - dx_{C_j^-} \wedge dy_{C_j^-} \wedge dx_{D_j^-} \wedge dy_{D_j^-}\,,$$
where $x$ and $y$ denote the real and imaginary parts, thus giving a factor $(-1)^{\wh g-1} = {\rm i}^{2\wh g-2}$ for $j = 2, \ldots, \wh g$. An extra factor of ${\rm i}$ similarly comes from the conversion of $dz_{d_1} \wedge d\bar{z}_{d_1}$ when setting $z_{c_1} = 1$.  Finally the extra factor $2$ is due to the fact that the hyperplane defined by $z_{c_1} = 1$ has distance $\sqrt{2}$ to the origin while the hyperplane defined by $z_{C_1^-} = 1$ has distance $1$ to the origin, hence the cone over the former has volume equal to $2$ times that of the latter, because the circle perimeter of radius $|z_{c_1}| = 1$ is also multiplied by $\sqrt{2}$ compared to that of radius $|z_{C_1^-}| = 1$.
\par
On the other hand, at the point $\sigma = (1, z_{d_1},0,\ldots,0)$ we have 
\bes
\omega_h \= \frac{1}{2\pi {\rm i}} \Biggl(\frac{\sum_{j = 2}^{\wt g} \chi_j
(dz_{c_j}\wedge d\bar{z}_{c_j} -  dz_{d_j}\wedge d\bar{z}_{d_j})}{h(\sigma)}
- \frac{\chi_1^2 \,dz_{d_1}\wedge d\bar{z}_{d_1}}{h(\sigma)^2} \Biggr)
\ees
and hence  
\begin{eqnarray*}
\nu_{h} & = & \frac{(-1)^{\wt{g}}(2\wt{g}-1)!}
{(2\pi {\rm i})^{2\wt{g}-1} h(\sigma)^{2\wt{g}}}
\Bigl(\chi_1^2 dz_{d_1}\wedge d\bar{z}_{d_1} \wedge \prod_{j =2}^{\wt g}
(\chi_j^2 dz_{c_j}\wedge d\bar{z}_{c_j} \wedge dz_{d_j}\wedge d\bar{z}_{d_j})
\Bigr)\, \\
& = &  \frac{(-1)^{\wt{g}}(2\wt{g}-1)!}
{2^{2g}(2\pi {\rm i})^{2\wt{g}-1} h(\sigma)^{2\wt{g}}}
\Bigl(dz_{d_1}\wedge d\bar{z}_{d_1} \wedge \prod_{j =2}^{\wt g}
(dz_{c_j}\wedge d\bar{z}_{c_j} \wedge dz_{d_j}\wedge d\bar{z}_{d_j})
\Bigr)\,,
\end{eqnarray*}
where we used the definition of $\chi_j$ in the last step. Comparing the above expressions of 
$\nu_{MV}$ and $\nu_{h}$ thus implies the desired identity.  
\end{proof}
\par
\begin{proof}[Proof of Theorem~\ref{thm:conjtrue}]
First note that by \cite{CoMoZa} the metric~$h$ on $\cOO_{\PP \wh{\cHH}}(-1)$
is good in the sense of Mumford \cite{mumford77}, and thus $\omega_h$ represents the first
Chern class of $\cOO_{\PP \wh{\cHH}}(-1)$ (after extending to the boundary). Moreover note that
$\zeta = 2\xi = c_1(\cOO_{\PP\wh \cHH}(2))$, which follows from the relation $q = \omega^2$ on the double cover. 
 Finally taking into account the order-$2$
stacky structure of $\PP\wh \cHH$ by the involution of the double cover, we conclude that  
\begin{eqnarray*}
\int_{\PP\barcQ[{g,s}](\nu)} \zeta^{2g-3+s}
& = & 2^{2g-2+s}
\int_{\PP \wh{\cHH}} \xi^{2g-3+s} \\
& = &  2^{2g-2+s} \, (-1)^{2g-3+s} \, \nu_h (\PP \wh\cHH) \\
& = & \frac{(2g-3+s)!}{2^{3-s}(2\pi {\rm i})^{2g-2+s}} \, 
\vol (\cQQ_{g,s}(\nu))\,,
\end{eqnarray*}
where we used Lemma~\ref{le:VFcomp} and the fact that $s$ is even in the last step. This thus implies the desired formula.
\end{proof}

%% file: sec_principal.tex
\section{Masur-Veech volumes of the principal strata}
\label{sec:principal}

In this section we evaluate the expression
given in Theorem~\ref{thm:conjtrue} for the principal
strata, and prove Theorem~\ref{thm:PSINT}, Theorem~\ref{thm:the7conj54} and Corollary~\ref{cor:volume-g=1} accordingly. For simplicity we abbreviate
the (fully labeled) principal stratum of quadratic differentials with $4g-4+n$ simple zeros and $n$ simple poles as $\cQQ_{g,N} = \cQQ_{g,N}(1^{4g-4+n},-1^n)$, where $N = 4g-4+2n$. Note the difference between $\cQQ_{g,N}$ and $\cQQ_{g,n}$, as the latter is the quadratic Hodge bundle over $\cMM_{g,n}$ with only the $n$ simple poles labeled.  In particular, their volumes differ by a factor $(4g-4+n)!$ due to the labeling of the $4g-4+n$ simple zeros.  

\subsection{The Segre class of $\BQ_{g,n}$} \label{sec:Segre}

Recall that for a projective bundle $p \colon \PP\cEE \to M$ associated to a vector bundle $\cEE$ of rank $r$, 
the $p$-pushforwards of $c_1(\cOO_{\PP\cEE}(1))^k$ are the Segre classes $s_{k-r+1}$ of~$\cEE$.
\par 
Here we consider the quadratic Hodge bundle $\BQ_{g,n}$ (extended over the Deligne-Mumford compactification) 
given by $\BQ_{g,n} = f_{*}(\omega_{f} (\sum_{i=1}^n \sigma_i))$, 
where $f\colon \cXX \to \BM_{g,n}$ is the universal curve, $\omega_{f}$ is the relative dualizing line bundle,  
and $\sigma_1, \ldots, \sigma_n$ are the sections of the $n$ marked points.  We denote by $\PP \BQ_{g,n}$ the associated projective bundle.  
Since the rank of $\BQ_{g,n}$ is $3g-3+n$ which equals $\dim \BM_{g,n}$, the top self-intersection number of the $\cOO(1)$-class $\zeta$ in Equation~\eqref{eq:odd} in this case thus corresponds to the degree of the top Segre class $s_{3g-3+n}(\BQ_{g,n})$, hence we can use it to compute  
the volumes of the principal strata.  
\par 
\begin{Prop}
\label{prop:volume-segre}
For the principal strata we have 
$$ \vol (\cQQ_{g,N}) = \frac{2^{2g+1} (\pi {\rm i})^{6g-6+2n} (4g-4+n)!}{(6g-7+2n)!} \int_{\BM_{g,n}} s_{3g-3+n} (\BQ_{g,n})\,.$$
\end{Prop}
\par 
\begin{proof}
The claim follows from Theorem~\ref{thm:conjtrue} for the special case $\nu = (1^{4g-4+n},-1^n)$ and $s = N = 4g-4+2n$, with an additional factor $(4g-4+n)!$ multiplied to the right-hand side, because in Equation~\eqref{eq:odd} all zeros are labeled while for the quadratic Hodge bundle $\cQQ_{g,n}$ over $\cMM_{g,n}$ we do not label the $4g-4+n$ simple zeros.  
\end{proof}
\par 
In order to evaluate the degree of $s_{3g-3+n}(\BQ_{g,n})$, we first compute the total Segre class of $\BQ_{g,n}$. The beginning step is a standard calculation 
originally due to \cite{mumford83} for the (Abelian) Hodge bundle and extended more generally in~\cite{Chiodo}, which we will further simplify. 
\par
Recall the definition
of the Bernoulli polynomials $B_n(x)$ by the expansion 
$$ \frac{te^{tx}}{e^t-1} \= \sum_{n \geq 0} \frac{B_n(x)}{n!} t^n\,. $$
They satisfy the properties that $B_n(x) = (-1)^n B_n(1-x)$ and that 
$$B_n(x+y) = \sum_{m=0}^n {n \choose m} B_m(x) y^{n-m}\,. $$
We also denote by $\kappa_d = f_{*}(c_1(\omega_f(\sum_{i=1}^n \sigma_i))^{d+1})$ for the $\kappa$-classes on $\BM_{g,n}$.  
\par
\begin{Lemma}[{\cite[Theorem~1.1.1]{Chiodo}} for the case $s=2$, $r=1$ and $m_i=1$] 
\label{lm:chern}
The Chern character of the quadratic Hodge bundle~$\BQ_{g,n}$ is given by
\begin{eqnarray*}
\ch (\BQ_{g,n}) & =  & 3g-3 +n +  \sum_{d\geq 1}\bigg( \frac{B_{d+1}(2)}{(d+1)!} \kappa_d - \sum_{i=1}^n \frac{B_{d+1}(1)}{(d+1)!}\psi_i^d \\
 & + & \frac{1}{2} \frac{B_{d+1}(1)}{(d+1)!}\sum_{i+j =d-1} i_{\irr *}(-\psi_1)^{i} \psi_2^j \\
 &+& \frac{1}{2} \frac{B_{d+1}(1)}{(d+1)!} \sum_{h=0}^{g}\sum_{S\subseteq [[1,n]]}\sum_{i+j =d-1} i_{h,S*} (-\psi_1)^{i} \psi_2^j\bigg)\,,   
\end{eqnarray*}
where the sum is constrained to $|S| \geq 2$ if $h = 0$ and to
$|S| \leq n-2$ if $h = g$.  In particular,  $i_{\irr}\colon \BM_{g-1, n+2}\to
\Delta_{\irr}$ has degree $2$ and the other $i_{h, S}\colon \BM_{h, S}\times
\BM_{g-h, S^c}\to \Delta_{h, S}$ are repeated twice in the sum (hence explaining the
factor $1/2$). 
\end{Lemma}
\par
Recall that for a vector bundle $\cEE$, the total Segre class
is given in terms of the coefficients of the Chern character by
\ba 
\label{eq:ch-segre}
 s(E; t) \= c(E;t)^{-1} \= \exp \Bigl(  \sum_{s \geq 1} (-1)^s (s-1)!
\ch_s (E) t^s  \Bigr)\,, 
\ea
where $t$ is the grading parameter. Hence in principle one can plug Lemma~\ref{lm:chern} in the above to compute the Segre classes of $\BQ_{g,n}$. 
 In order to further simplify the
expression obtained this way,
we introduce the notation $\kappa_{(d)}=\pi_*(\psi_{n+1}^2\ldots \psi_{n+d}^2)$,
where $\pi$ is the map $\overline{\mathcal{M}}_{g,n+d}\to
\overline{\mathcal{M}}_{g,n}$ induced by forgetting the last $d$ marked points.
\par
\begin{Lemma} 
\label{le:segreBQ}
The total Segre class of $\BQ_{g,n}$ can be expressed as   
$$
s(\BQ_{g,n}) \= \bigl(1-\lambda_1+\cdots + (-1)^g \lambda_g\bigr)
\left(1- \frac{\kappa_{(1)}}{1!}
+ \frac{\kappa_{(2)}}{2!}- \frac{\kappa_{(3)}}{3!}+ \cdots\right)\,.
$$
In particular, the top Segre class of $\BQ_{g,n}$ is equal to   
$$ 
s_{3g-3+n}(\BQ_{g,n}) \=  (-1)^{3g-3+n} \sum_{i=0}^g \frac{\kappa_{(2g-3+n+i)} \lambda_{g-i}}{(2g-3+n+i)!}\,.
$$
\end{Lemma}
\par
\begin{proof} 
Denote by $\BE_{g,n}$ the (Abelian) Hodge bundle over $\BM_{g,n}$. 
Comparing the expression of $\ch (\BQ_{g,n})$ in Lemma~\ref{lm:chern} with that of $\ch (\BE_{g,n})$ in~\cite[Equation (1)]{Chiodo}, we conclude that 
\ba
\label{eq:E-Q}
 \ch (\BQ_{g,n}) & \= \ch (\BE_{g,n}) - 1 + \sum_{d\geq 0} \Big( \frac{B_{d+1}(2) - B_{d+1}(1)}{(d+1)!} \kappa_d \Big) \\
 & \= \ch (\BE_{g,n}) - 1 + \sum_{d\geq 0}\frac{1}{d!} \kappa_d\,. 
\ea
\par
To further simplify this expression we recall some $\kappa$-class computation from
\cite[Section~2]{PixtonThesis}. For any permutation $\sigma\in S_m$
with cycle decomposition $\sigma = \sigma_1 \cdots \sigma_\ell$
and integers $\alpha_1,\ldots,\alpha_m$, we define $\kappa^\alpha_\sigma
= \prod_{j} \kappa_{|\sigma_j|}$ where~$|\sigma_j| = \sum_{i \in \sigma_j} \alpha_i$
is the sum of~$\alpha_i$ corresponding to the cycle. Then there is the
pushforward formula
\bes
\pi_*(\psi_{n+1}^{\alpha_1+1} \ldots \psi_{n+m}^{\alpha_m+1})
\= \overbrace{\kappa_{\alpha_1}\cdots \kappa_{\alpha_m}} \,\coloneqq \,
\sum_{\sigma \in S_m} \kappa_\sigma^{\alpha}\, 
\ees
(see~\cite[Section 2.3]{PixtonThesis}). 
In particular, in this notation $\kappa_{(j)} = \overbrace{\kappa_{1}\cdots
\kappa_{1}}$ with $j$ factors. On the other hand, we define the linear map
\bes
\Bigl\{\sum_i c_i t^i \Bigr\}_\kappa \= \sum_i c_i \kappa_i t^i 
\ees
on power series. Then~\cite[Lemma~2.3]{PixtonThesis}
shows that 
\be \label{eq:pixton23}
\exp ( - \{ \log (1-X) \}_\kappa) \= \overbrace{\exp(\{ X\}_\kappa)}\, 
\ee
for any $X \in \QQ[t]$.  Combining~\eqref{eq:ch-segre} with~\eqref{eq:E-Q} thus implies that 
\begin{eqnarray*}
s(\BQ_{g,n}) & \= & s(\BE_{g,n}) \, \exp \Big( \sum_{d \geq 1} \frac{(-1)^d}{d} \kappa_d   \Big)  \\
& \= & c(\BE_{g,n}^{\ast}) \, \exp\Big(-\kappa_1 + \frac{\kappa_2}{2} -  \frac{\kappa_3}{3} + \cdots \Big) \\
& \= &  \big(1 - \lambda_1 + \lambda_2 - \cdots + (-1)^g \lambda_g\big) \left(1- \frac{\kappa_{(1)}}{1!}+ \frac{\kappa_{(2)}}{2!}- \frac{\kappa_{(3)}}{3!}+ \cdots\right)\,, 
\end{eqnarray*}
where we used~\eqref{eq:pixton23} for $X(t) = -t$ and the fact that $s(\BE_{g,n}) = c(\BE_{g,n})^{-1} = c(\BE_{g,n}^{\ast})$ (see~\cite[Equation (5.4)]{mumford83}). 
\par 
The expression of $s_{3g-3+n}(\BQ_{g,n})$ follows from taking the terms of codimension $3g-3+n$ in the expansion of the above product.  
\end{proof}

\subsection{Linear Hodge integrals}
 
In this section we use the previous results to prove Theorem~\ref{thm:PSINT}. To simplify notation, for $n\geq 0$ and $0\leq i \leq g$ we define the following linear Hodge integrals 
\begin{eqnarray}
\label{eq:def-kappa-n}
\varkappa(g,i)_n & \coloneqq &  \int_{\BM_{g,n}} \kappa_{(2g-3+n+i)} \lambda_{g-i} \\ \nonumber
& \= & \int_{\BM_{g, 2g-3+2n+i}} \psi_{n+1}^2 \cdots \psi_{2g-3+2n+i}^2 \lambda_{g-i}\,,
\end{eqnarray}
where the second equality follows from the definition of $\kappa_{(d)}$ and the projection formula.  
For $n=0$ we also denote $\varkappa(g,i)_0=\varkappa(g,i)$, i.e. 
\begin{eqnarray}
\label{eq:def-kappa-0}
\varkappa(g,i)  & \coloneqq &  \int_{\BM_{g, 2g-3+i}} \psi_{1}^2 \cdots \psi_{2g-3+i}^2 \lambda_{g-i}\,. 
\end{eqnarray}
We show that $\varkappa(g,i)_n$ can be expressed in terms of $\varkappa(g,i)$. 
\par
 \begin{Lemma}
 \label{lm:n-to-0}
 For all $g$, $i$ and $n$ we have 
 $$
\varkappa(g,i)_n \= \varkappa(g,i) \, \frac{(2g-3+n+i)!}{(2g-3+i)!}\, \frac{(4g-7+2n+i)!!}{(4g-7+i)!!}\,.
 $$
 \end{Lemma}
 \par
 \begin{proof} The proof relies on the string and dilation equations \cite[(2.41) and (2.45)]{witten}.
We choose a triple $(g,i,n)$ such that $n>0$ and denote by $d=2g-3+n+i$.
Then we have 
\begin{eqnarray}
\label{eq:kappa-string}
 \varkappa(g,i)_n &=& \int_{\oM_{g,n+d}} \psi_{n+1}^2 \ldots \psi_{n+d}^2 \lambda_{g-i}\\ \nonumber 
 &=&  \sum_{j=1}^d \int_{\oM_{g,n-1+d}} \psi_{n}^2 \ldots \psi_{n-1+j} \ldots \psi_{n+d-1}^2 \lambda_{g-i}\\  \nonumber
 &=& d   \int_{\oM_{g,n-1+d}} \psi_{n}^2 \ldots \psi_{n+d-2}^2 \psi_{n+d-1} \lambda_{g-i}\\  \nonumber
 &=& d  (2g-4+d+n)  \int_{\oM_{g,n-2+d}} \psi_{n}^2 \ldots \psi_{n+d-2}^2  \lambda_{g-i} \\  \nonumber
 &=& d  (2d-i-1) \varkappa(g,i)_{n-1}\,,
 \end{eqnarray}
where from the first line to the second we used the string equation applied
to the first marked point and from the third line to the fourth we used the dilation
equation for the last marked point. The claim thus follows by induction on~$n$.
 \end{proof}
\par 
\begin{proof}[Proof of Theorem~\ref{thm:PSINT}]
Proposition~\ref{prop:volume-segre} and Lemma~\ref{le:segreBQ} imply that 
$$ 
\vol (\cQQ_{g,N}) =  \frac{2^{2g+1} (\pi {\rm i})^{6g-6+2n} (4g-4+n)!}{(6g-7+2n)!} (-1)^{3g-3+n}\sum_{i=0}^g \frac{\varkappa(g,i)_n}{(2g-3+n+i)!}\,,
$$ 
which after simplification gives the first equality in~\eqref{eq:principal}. The second equality in~\eqref{eq:principal} thus follows from 
Lemma~\ref{lm:n-to-0}. 
\end{proof}

\subsection{Volumes of the principal strata for fixed $g$
  and varying $n$}
  \label{subsec:volume-large-n}

In this section we consider numerical properties of $\vol(\cQQ_{g,N})$ when $g$ is fixed and $n$ varies (i.e. $N = 4g-4+2n$ varies), especially as $n$ tends to infinity.  
\par
\begin{proof}[Proof of Theorem~\ref{thm:the7conj54}, Equation~\eqref{eq:volume-fix-g}]
We renormalize the volume to be 
\begin{eqnarray}
\label{eq:v(g,n)}
v(g,n)  \=  \frac{\vol(\cQQ_{g,N})}{2^n  \pi^{6g-6+2n}}
\, \frac{(6g-7+2n)!}{(2g-3+n)!(4g-4+n)!} \,.
\end{eqnarray}
Then Theorem~\ref{thm:PSINT} implies that 
\begin{eqnarray*}
2^{2-4g} v(g,n) 
 \=  \frac{1}{(4g-6+2n)!!} \sum_{i=0}^g \frac{\varkappa(g,i)}{(2g-3+i)!} \, \frac{(4g-7+2n+i)!!}{(4g-7+i)!!}\,.
\end{eqnarray*}
Denote by 
$$\varkappa(g,i)' \= \frac{\varkappa(g,i)}{(2g-3+i)!(4g-7+i)!!}$$ 
which is independent of $n$. Setting 
$a_{g,n} = 4g-6+2n$ below, we conclude that 
\bas \
 & \phantom{\=} 2^{2-4g} v(g,n) \,\=\,  \sum_{i=0}^g \varkappa(g,i)'
\frac{(a_{g,n}+i-1)!!}{a_{g,n}!!}\\ 
& \=  \sum_{i=0}^{\lfloor (g-1)/2\rfloor}  \varkappa(g,2i +1)'
\frac{(a_{g,n}+ 2i)!!}{a_{g,n} !!} +   \sum_{i=0}^{\lfloor g/2\rfloor}
\varkappa(g, 2i)' \frac{(a_{g,n}-1 + 2i)!!}{a_{g,n}!!}\\
& \=  \sum_{i=0}^{\lfloor (g-1)/2\rfloor}  \varkappa(g, 2i+1)'\, \frac{(a_{g,n}+2i)!!}
{a_{g,n}!!} +  \gamma_{2g-3+n} \sum_{i=0}^{\lfloor g/2\rfloor}
\varkappa(g, 2i)'\,\frac{(a_{g,n}-1+2i)!!}{(a_{g,n} -1)!!}\,,
\eas
where in the last step we used $\gamma_k = \frac{1}{4^k} {2k \choose k}$. Therefore, we obtain that 
$$ p_g(n) \= 2^{4g-2}\sum_{i=0}^{\lfloor (g-1)/2\rfloor}  \varkappa(g, 2i+1)'\frac{(4g-6+2n+2i)!!}{(4g-6+2n)!!} $$
and
$$ q_g(n) \= 2^{4g-2}\sum_{i=0}^{\lfloor g/2\rfloor} \varkappa(g, 2i)'\frac{(4g-7+2n+2i)!!}{(4g-7+2n)!!} $$
as the polynomials whose existence we claimed to exist. 
\par 
Note that the leading coefficients of $p_g$ and $q_g$ are given by $\int_{\BM_{g,3g-3}}\psi_1^2 \cdots \psi_{3g-3}^2$ and $\int_{\BM_{g,3g-4}}\psi_1^2 \cdots \psi_{3g-4}^2\lambda_1$. It is well-known that the divisor class $\lambda_1$ is nef and effective and that the divisor classes $\psi_i$ are nef (see e.g.~\cite[Chapter 6, D]{harrismorrison}). Bigness of $\psi_i$ also holds, since the top self-intersection of $\psi_i$ on $\BM_{g,n}$ is positive (see~\cite[Equation (28)]{FPHodge1}). Then the non-vanishing of the lead terms of $p_g$ and $q_g$ follows from Lemma~\ref{le:positivity} below.
\end{proof} 
\par
\begin{Lemma} \label{le:positivity}
Suppose $D_1, \ldots, D_n$ are divisor classes in an $n$-dimensional projective
variety~$X$ such that~$D_1$ is nef and effective and such that
$D_2, \ldots, D_n$
are big and nef. Then the intersection number $\int_X D_1\cdots D_n$ is positive.  
\end{Lemma}
\par
\begin{proof}
Since $D_n$ is big, we can write $D_n = A_n + E_n$, where $A_n$ is ample
and $E_n$ is effective.  Then $\int_X D_1\cdots D_{n-1}  E_n \geq 0$,
because a nef divisor class is a limit of ample divisor classes.  Hence it suffices
to show that $\int_X D_1\cdots D_{n-1} A_n > 0$.  Repeating the
argument inductively, it eventually reduces to showing that
$\int_X D_1 A_2 \cdots A_n > 0$, where the $A_i$ are all ample,
and this holds since $D_1$ is effective.
\end{proof}
\par 
 We can extract the leading term of the volume from the above expression. Recall that $2^{6g-7} m_g$ is the top coefficient 
 of $q_g$ if $g$ is even or the top coefficient of $p_g$ if $g$ is odd, and that $\epsilon (g) = 0$ or $1$ determined by the parity of $g$.
 \par
 \begin{proof}[Proof of Theorem~\ref{thm:the7conj54}, Equation~\eqref{eq:volume-large-n}]
Note that $\gamma_k \sim (\pi k)^{-1/2}$ as $k\to \infty$. Then
as $n\to \infty$, the dominant term of the two in the sum  
$v(g,n) = p_g(n) + \gamma_{2g-3+n}q_g(n)$ is $\gamma_{2g-3+n}q_g(n)$ when~$g$
is even and is $p_g(n)$ when~$g$ is odd. 
 Hence the leading term of $v(g,n)$ (as a function of $n$) is $\pi^{-1/2} n^{(g-1)/2} 2^{6g-7} m_g$ when $g$ is even and is $n^{(g-1)/2} 2^{6g-7} m_g$ when $g$ is odd. Altogether it implies that for large $n$ 
 \ba
 \label{eq:v(g,n)-asy}
 v(g,n) \sim   2^{6g-7} \pi^{(\epsilon (g)-1)/2} m_g n^{(g-1)/2}\,.\\
 \ea
 For the leading term of ${6g-7 + 2n \choose 2g-3+n}$, recall Stirling's approximation $k! \sim \sqrt{2\pi k} (k/e)^k$. It implies that for large $n$ 
$$ {6g-7 + 2n \choose 2g-3+n}\sim 2^{2n+6g-7} \frac{1}{\sqrt{\pi n}}\,. $$
The claim on $\vol(\cQQ_{g,N})$ thus follows from~\eqref{eq:v(g,n)-asy} and 
the conversion of $\vol(\cQQ_{g,N})$ to $v(g,n)$ in~\eqref{eq:v(g,n)}.   
\end{proof}
\par 
From the above proof, we see that the coefficient $m_g$ is an explicit rescaling of the intersection number $\int_{\BM_{g,3g-3}} \psi_1^2 \cdots \psi_{3g-3}^2$. The generating series of such intersection numbers satisfies the Painlev\'e equation~I, and hence they can be computed efficiently. We refer to~\cite[Section~6]{ItZuKont}, \cite[Section~4.2]{zvonkine}, \cite[Section 4]{LMXHodge}, \cite{DYZHur}, and~\cite{YZZVolume} for related discussions on this topic.  

 \subsection{Volumes of the principal strata in genus one}
\label{subsec:volume-g=1}
 
In this section we prove Corollary~\ref{cor:volume-g=1}. By Lemma~\ref{le:segreBQ} we have 
$$ \int_{\BM_{1,n}} s_{n}(\BQ_{1,n}) = (-1)^n  \int_{\BM_{1,n}} \left(   \frac{\kappa_{(n)}}{n!}
+ \frac{\lambda_1 \kappa_{(n-1)}}{(n-1)!}\right)\,. $$
 By the $\lambda_g$-theorem (\cite{FPHodge1}), we can evaluate the Hodge integral 
 \ba
 \label{eq:segre-next-g=1}
  \int_{\BM_{1,n}}\lambda_1 \kappa_{(n-1)} & \=  \int_{\BM_{1,2n-1}} \psi_{n+1}^{2}\cdots \psi_{2n-1}^2 \lambda_1 \\ 
  & \=  {2n-2 \choose 0, \ldots, 0, 2, \ldots, 2} \frac{1}{24}
  \=  \frac{(n-1)!(2n-3)!!}{24}\,.
  \ea
Moreover, using the string and dilation equations we have
\bas
   \int_{\BM_{1,n}} \kappa_{(n)} & \=  \int_{\BM_{1,2n}} \psi_{n+1}^{2}\cdots \psi_{2n}^2 \\
   & \=  n \int_{\BM_{1,2n-1}} \psi_{n}^2 \cdots \psi_{2n-2}^2 \psi_{2n-1} \\
   & \=  n(2n-2)  \int_{\BM_{1,2n-2}}  \psi_{n}^2 \cdots \psi_{2n-2}^2
   \= \cdots \=  \frac{1}{24} n! (2n-2)!!\,.
\eas
It follows that 
\ba
\label{eq:segre-g=1}
\int_{\BM_{1,n}} s_{n}(\BQ_{1,n})  = (-1)^n  \frac{1}{24} \big((2n-3)!!
+ (2n-2)!!\big)\,.
\ea
By Proposition~\ref{prop:volume-segre} we thus conclude that 
  \begin{eqnarray*}
  \vol(\cQQ_{1,2n}(1^n, -1^n)) & \= & \frac{2^{3}(\pi {\rm i})^{2n} n!}{(2n-1)!}\, s_{n}(\BQ_{1,n})  \\
  & \= & \pi^{2n} \frac{n!}{3 (2n-1)!} \big((2n-3)!! + (2n-2)!!\big)\,,
\end{eqnarray*}
which proves Equation~\eqref{eq:volume-g=1}. In particular, it confirms~\cite[Conjecture 5.4]{TRMVV} for the case $g=1$ with $p_1(n) = q_1(n) = 1/6$.

%% file: sec_svLyap.tex
\section{Siegel-Veech constants and Lyapunov exponents}
\label{sec:svL}

In this section we study area Siegel-Veech constants $c_{\area}$ and sums of (involution-invariant) Lyapunov exponents $L^{+}$ for the principal strata of quadratic differentials.  In particular, we will prove Theorem~\ref{thm:sv-L}, Corollary~\ref{cor:sv-L}, and Corollary~\ref{cor:sv-L-g=1}. In addition, we will give conjectural formulas to compute area Siegel-Veech constants and sums of Lyapunov exponents as intersection numbers for arbitrary affine invariant submanifolds. 

\subsection{$c_{\area}$ and $L^{+}$ as intersection numbers}
\label{subsec:sv-L-intersection}

Recall that $\barcQ[g,n]$ is the quadratic Hodge bundle over $\BM_{g,n}$ with only $n$ marked points, and $\delta$ is the divisor class of 
the total boundary of $\PP\barcQ[g,n]$. We first prove the formula for $c_{\area}$. 
\par 
\begin{proof}[Proof of Theorem~\ref{thm:sv-L}, Equation~\eqref{eq:sv-principal}]
For simplicity we denote by $s_{g,n}$ the degree of the top Segre class of $\barcQ[g,n]$. Note that the denominator of the right-hand side of~\eqref{eq:sv-principal} corresponds to $s_{g,n}$. By \cite[Section 4.2, Corollary 1]{GouSV} (see also \cite[Theorem~4.1]{TRMVV}) and Proposition~\ref{prop:volume-segre}, the desired formula is equivalent to the following equality of top Segre classes 
\ba
\label{eq:sv-segre}
 \int_{\PP\barcQ[g,n]}\zeta^{6g-8+2n}\delta & \= \frac{1}{2} s_{g-1,n+2} + \frac{1}{2} \sum_{g_1+g_2 = n \atop n_1+n_2 = n} \frac{n!}{n_1! n_2!} 
 s_{g_1,n_1+1}s_{g_2,n_2+1}\,,
\ea
where the sum ranges over admissible pairs $(g_i, n_i)$ as constrained in Lemma~\ref{lm:chern}. Since the total boundary of $\BM_{g,n}$ is the union of irreducible boundary divisors whose types correspond to the summands on the right-hand side of \eqref{eq:sv-segre}, the equality follows from the structure of $\barcQ[g,n]$ restricted to each of the boundary divisors.  
\par 
More precisely, consider the morphism $i_{\irr}\colon \BM_{g-1,n+2} \to \Delta_{\irr}$ induced by normalizing a non-separating node $r$ of a pointed stable curve $(X, p_1,\ldots, p_n)$. Let $X'$ be the normalization of $X$ at $r$ and denote by $r_1, r_2$ the (labeled) preimages of $r$ in $X'$. The fiber of $\barcQ[g,n]$ over $(X, p_1,\ldots, p_n)$ parameterizes quadratic differentials $q \in H^0(\omega_{X'}^{\otimes 2} (p_1 + \cdots + p_n + 2r_1 + 2r_2))$ such that the $2$-residues of $q$ at $r_1$ and $r_2$ are equal (see~\cite[Section 3.1]{kdiff} for $k$-residues in general). From this viewpoint we can regard $\barcQ[g-1,n+2]$ as a sub-bundle of $i_{\irr}^{*}\barcQ[g,n]$ whose fiber over $(X', p_1,\ldots, p_n, r_1, r_2)$ is $H^0(\omega_{X'}^{\otimes 2} (p_1 + \cdots + p_n + r_1 + r_2))$, containing quadratic differentials with zero $2$-residues at $r_1$ and $r_2$. We then have the exact sequence 
$$ 0 \to \barcQ[g-1,n+2] \to i_{\irr}^{*}\barcQ[g,n] \to \mathbb{C} \to 0\,, $$ 
where the last map to the trivial line bundle is induced by taking the value of the $2$-residue at $r_1$. It implies that 
$$s_{g-1,n+2} = \int_{\BM_{g-1,n+2}} s_{3g-4+n}(i_{\irr}^{*}\barcQ[g,n]) = 2 \int_{\PP\barcQ[g,n]}\zeta^{6g-8+2n}\delta_{\irr}\,,$$ 
where the factor $2$ is due to the labeling of $r_1$ and $r_2$, i.e. due to $\deg i_{\irr} = 2$. This explains the term $\frac{1}{2} s_{g-1,n+2}$ on the right-hand side of~\eqref{eq:sv-segre}. 
\par 
Similarly consider the morphism $i_{g_1,S_1}\colon \BM_{g_1, n_1}\times \BM_{g_2, n_2} \to \Delta_{g_1, S_1}$ induced by normalizing a separating node $r$ of a pointed stable curve $(X, p_1,\ldots, p_n)$, where $X$ consists of two components $X_1$ and $X_2$ separated by $q$,   
$g_1 + g_2 = g$, $n_1 + n_2 = n$, and $X_i$ contains the $n_i$ marked points in the subset $S_i \subseteq \{ p_1, \ldots, p_{n}\}$ for $i = 1, 2$. 
Using the same argument as in the previous case, we obtain the exact sequence 
$$ 0 \to  \barcQ[g_1,n_1+1] \oplus \barcQ[g_2,n_2+1] \to i_{g_1, S_1}^{*} \barcQ[g,n] \to \mathbb{C} \to 0\,. $$
It implies that 
$$ \int_{\PP\barcQ[g,n]}\zeta^{6g-8+2n}\delta_{g_1, S_1} = \int_{\BM_{g_1, n_1}\times \BM_{g_2, n_2}} s_{3g-4+n} ( i_{g_1, S_1}^{*} \barcQ[g,n]) =  s_{g_1,n_1+1}s_{g_2,n_2+1}\,.$$
This explains the remaining summation on the right-hand side of~\eqref{eq:sv-segre}, where the factors $\frac{1}{2}$ and $\frac{n!}{n_1! n_2!}$ are due to the ordering of $n_1, n_2$ and choosing $n_1$ marked points out of the total $n$ marked points respectively.  
\end{proof}
\par
Next we prove the formula for $L^+$. 
\par 
\begin{proof}[Proof of Theorem~\ref{thm:sv-L}, Equation~\eqref{eq:L-principal}]
Denote by $\kappa_{g,n} = \frac{1}{18}(5g-5-n)$. Then by~\cite[Theorem 2, (2.3)]{ekz} we have the following relation 
\ba
\label{eq:ekz}
 L^+(1^{4g-4+n}, -1^n)  \= \kappa_{g,n} + \frac{\pi^2}{3} c_{\area}(1^{4g-4+n}, -1^n)\,.
 \ea
For a family of nodal curves $f\colon \mathcal{X} \to B$, there is a relation of divisor classes $12 \lambda_1 - \delta  = f_{*}(c_1 (\omega_{f})^2)$ (see~\cite[p. 306]{mumford83}). Combining with~\eqref{eq:sv-principal}, the desired formula is equivalent to the following equality 
\ba
\label{eq:kappa}
\int_{\PP\barcQ[g,n]}\zeta^{6g-8+2n} f_{*}(c_1 (\omega_{f})^2) & \= - 6\kappa_{g,n} \int_{\PP\barcQ[g,n]}\zeta^{6g-7+2n}\,, 
\ea  
where $f\colon \mathcal{X} \to \PP\barcQ[g,n]$ is the universal curve.  
\par 
Let $\PP\barcQ[g,n](2)$ be the closure of the locus in $\PP\barcQ[g,n]$ parameterizing quadratic differentials with a double zero. A general differential parameterized in $\PP\barcQ[g,n](2)$ has zero type $(2, 1^{4g-6+n}, -1^n)$. Let $\PP\barcQ[g,n](0)$ be the closure of the locus in $\PP\barcQ[g,n]$ parameterizing quadratic differentials that are holomorphic at some marked point $p_i$.  A general differential parameterized in $\PP\barcQ[g,n](0)$ has zero type $(1^{4g-5+n}, -1^{n-1}, 0)$, where the entry $0$ indicates that $p_i$ becomes an ordinary point instead of a simple pole (still labeled). Denote by $B$ the complement of 
$\PP\barcQ[g,n](2)$ and $\PP\barcQ[g,n](0)$. Then differentials parameterized by $B$ have zero type exactly $(1^{4g-4+n}, -1^n)$. 
\par 
Denote by $Z_i$ and $P_j$ the $i$-th zero section and the $j$-th pole section in the universal curve $\mathcal{X}$ over $B$ respectively, which are pairwise disjoint. Then as in~\cite[Section 3.4]{ekz} we have the relation of divisor classes 
$$ f^{*} (-\zeta) = 2 c_1(\omega_{f}) - \sum_{i=1}^{4g-4+n} Z_i + \sum_{j=1}^n P_j \,.$$ 
Intersecting the relation with $Z_i$ and $P_j$ respectively and pushing forward by $f$, we conclude that 
$ f_{*}(Z_i^2) = \zeta / 3$ and $f_{*}(P_j^2) = \zeta $ for all $i$ and $j$. Intersecting the relation with $c_1(\omega_f)$ and pushing forward by $f$, we  obtain that $f_{*}(c_1 (\omega_{f})^2) = - 6\kappa_{g,n} \zeta$ in $B$, which implies that the same relation holds in $\PP\barcQ[g,n]$ modulo a divisor class supported on the union of $\PP\barcQ[g,n](2)$ and $\PP\barcQ[g,n](0)$. Therefore, it suffices to prove that 
$\int_{\PP\barcQ[g,n](2)} \zeta^{6g-8+2n} = 0$ and $\int_{\PP\barcQ[g,n](0)} \zeta^{6g-8+2n} = 0$. 
\par 
Note that $\PP\barcQ[g,n](2)$ is the closure of $\PP\cQQ(2, 1^{4g-6+n}, -1^n)$ whose period coordinates are not all given by absolute periods due to the double zero. Recall in Section~\ref{sec:IntMV} that $-\zeta$ corresponds to the curvature form $\omega_h$ of the good metric on the tautological bundle $\cOO(-1)$. Then $\omega_h^{6g-8+2n} = 0$ on this stratum, since the exponent is one bigger than the rank of the subspace of absolute periods (after projectivization). It follows that $\int_{\PP\barcQ[g,n](2)} \zeta^{6g-8+2n} = 0$.  Similarly $\PP\barcQ[g,n](0)$ consists of closures of $\PP\cQQ(1^{4g-5+n}, -1^{n-1}, 0)$ whose period coordinates are not all given by absolute periods due to the labeled ordinary point (as a zero of order $0$).  Then the same argument implies that $\int_{\PP\barcQ[g,n](0)} \zeta^{6g-8+2n} = 0$. 
\end{proof}

\subsection{$c_{\area}$ and $L^{+}$ for fixed $g$ and varying $n$}
\label{subsec:sv-L-large-n}

In this section we consider numerical properties of $c_{\area}$ and $L^+$ for the principal strata when $g$ is fixed and $n$ varies, especially as $n$ tends to infinity. Since $c_{\area}$ and $L^+$ determine each other by the relation~\eqref{eq:ekz}, it suffices to evaluate $L^+$.  
\par
The denominator in the formula~\eqref{eq:L-principal} for $L^+$ corresponds to the top Segre class of $\barcQ[g,n]$, which has been computed via linear Hodge integrals of type~\eqref{eq:def-kappa-n} and~\eqref{eq:def-kappa-0}. To evaluate the numerator, we define similarly 
\ba
\label{eq:def-theta-n}
 \vartheta(g, i)_n & \coloneqq  \int_{\BM_{g,n}} \kappa_{(2g-4+n+i)}\lambda_{g-i}\lambda_1 \\ 
 & \=  \int_{\BM_{g, 2g-4+2n+i}} \psi_{n+1}^2 \cdots \psi_{2g-4+2n+i}^2 \lambda_{g-i}\lambda_1\,. 
\ea
For $n = 0$ we also denote $\vartheta (g, i)_0 = \vartheta(g, i)$, i.e. 
\ba
\label{eq:def-theta-0}
\vartheta(g, i) & \coloneqq & \int_{\BM_{g,2g-4+i}} \psi_1^2 \cdots \psi_{2g-4+i}^2 \lambda_{g-i}\lambda_1\,.
\ea
We show that $\vartheta(g,i)_n$ can be expressed in terms of $\vartheta(g,i)$. 
\par
\begin{Lemma}
\label{lm:theta-n-to-0}
For all $g$, $i$ and $n$ we have 
$$ \vartheta(g, i)_n \= \vartheta(g, i)\, \frac{(2g-4+n+i)!}{(2g-4+i)!} \, \frac{(4g-8+2n+i)!!}{(4g-8+i)!!}\,. $$
\end{Lemma}
\par
\begin{proof}
Let $d = 2g - 4 + n + i$. For any $n > 0$, the same method as in the proof of Lemma~\ref{lm:n-to-0} implies that 
\begin{eqnarray*}
 \vartheta(g, i)_n & = &  \int_{\BM_{g, n+d}} \psi_{n+1}^2 \cdots \psi_{n+d}^2 \lambda_{g-i}\lambda_1 \\
 & = & d (2d - i ) \vartheta(g, i)_{n-1}\,. 
 \end{eqnarray*}
 Then the claim follows by induction on $n$. 
\end{proof}
\par
\begin{Lemma}
\label{lm:segre-next}
The numerator of~\eqref{eq:L-principal} is equal to 
$$ \int_{\BM_{g,n}} s_{3g-4+n}(\BQ_{g,n}) \lambda_1 \= (-1)^{3g-4+n} \sum_{i=0}^g \frac{\vartheta(g,i)_n}{(2g-4+n+i)!}\,. $$
\end{Lemma}
\par
\begin{proof}
This follows from the expression of the total Segre class $s(\BQ_{g,n})$ in Lemma~\ref{le:segreBQ} and the definition of $\vartheta(g,i)_n$ in~\eqref{eq:def-theta-n}.
\end{proof}
\par 
Now we can verify the desired numerical properties of $L^+$ and $c_{\area}$ for fixed $g$ and varying $n$. 
\par
\begin{proof}[Proof of Corollary~\ref{cor:sv-L}]
We renormalize the numerator of~\eqref{eq:L-principal} to be 
\ba
\label{eq:u(g,n)}
u(g, n) \= \frac{(-1)^{3g-4+n} 2^{4g-2}}{(4g-8+2n)!!}  \int_{\BM_{g,n}} s_{3g-4+n}(\BQ_{g,n}) \lambda_1\,. 
\ea
Then Lemma~\ref{lm:theta-n-to-0} and Lemma~\ref{lm:segre-next} imply that 
\bas
2^{2-4g} u(g,n) &\=   \frac{1}{(4g-8+2n)!!}  \sum_{i=0}^g \frac{\vartheta(g,i)_n}{(2g-4+n+i)!} \\
 \ &\=  \frac{1}{(4g-8+2n)!!}  \sum_{i=0}^g \frac{\vartheta(g,i)}{(2g-4+i)!} \, \frac{(4g-8+2n+i)!!}{(4g-8+i)!!}\,.
\eas
Denote by 
$$ \vartheta(g,i)' \= \frac{\vartheta(g,i)}{(2g-4+i)! (4g-8+i)!!}\, $$
which is independent of $n$. Setting $b_{g,n} = 4g-8+2n$ below, we conclude that 
\begin{eqnarray*}
2^{2-4g} u(g,n)& = &  \sum_{i=0}^g \vartheta(g, i)' \frac{(b_{g,n}+i)!!}{b_{g,n}!!} \\
& = & \sum_{i=0}^{\lfloor (g-1)/2\rfloor}  \vartheta(g,2i +1)'\frac{(b_{g,n}+2i+1)!!}{b_{g,n}!!} +   \sum_{i=0}^{\lfloor g/2\rfloor} \vartheta(g, 2i)' \frac{(b_{g,n}+2i)!!}{b_{g,n}!!}\\
& = &(4g-6+2n)\gamma_{2g-3+n}\sum_{i=0}^{\lfloor (g-1)/2\rfloor}  \vartheta(g,2i +1)'\frac{(b_{g,n}+2i+1)!!}{(b_{g,n}+1)!!}\\
& + &  \sum_{i=0}^{\lfloor g/2\rfloor} \vartheta(g, 2i)' \frac{(b_{g,n}+2i)!!}{b_{g,n}!!}\,. 
\end{eqnarray*}
Therefore, we obtain that 
$$ r_g(n) \= 2^{4g-2} \sum_{i=0}^{\lfloor g/2\rfloor} \vartheta(g, 2i)' \frac{(4g-8+2n+2i)!!}{(4g-8+2n)!!}\, $$
and 
$$ s_g(n) \= 2^{4g-2} (4g-6+2n) \sum_{i=0}^{\lfloor (g-1)/2\rfloor}  \vartheta(g,2i +1)'\frac{(4g-7+2n+2i)!!}{(4g-7+2n)!!}\, $$
as the polynomials whose existence we claimed to exist. The non-vanishing of their leading coefficients follows from a similar argument using Lemma~\ref{le:positivity}. 
\par
Recall $v(g,n)$ defined in~\eqref{eq:v(g,n)} in the proof of Theorem~\ref{thm:the7conj54}. In summary, we have 
\bas
 v(g,n) &\=   \frac{(-1)^{3g-3+n} 2^{4g-2}}{(4g-6+2n)!!} \int_{\BM_{g,n}} s_{3g-3+n}(\BQ_{g,n}) \=  p_g(n) + \gamma_{2g-3+n}q_g(n)\,, \\
 u(g,n) & \=  \frac{(-1)^{3g-4+n}  2^{4g-2}}{(4g-8+2n)!!} \int_{\BM_{g,n}} s_{3g-4+n}(\BQ_{g,n}) \lambda_1 \=  r_g(n) + \gamma_{2g-3+n}s_g(n)\,. 
 \eas
The rational function expression of $L^+$ in~\eqref{eq:L-fixed-g-large-n} thus follows from~\eqref{eq:L-principal}.  
\par 
For the large $n$ asymptotic of $L^+$, recall that $2^{6g-7} n_g$ is the leading coefficient of $r_g(n)$ and $s_g(n)$ when $g$ is even and odd respectively. We claim that 
\ba
\label{eq:u(g,n)-asy}
u(g,n) \sim 2^{6g-7} \pi^{-\epsilon (g)/2} n_g n^{g/2}\, 
\ea
for fixed $g$ and $n\to \infty$. To see this, the dominant term of the two in the sum  
 $u(g,n) = r_g(n) + \gamma_{2g-3+n}s_g(n)$ is $r_g(n)$ when $g$ is even and $\gamma_{2g-3+n}s_g(n)$ when $g$ is odd. 
 Hence the leading term of $u(g,n)$ (as a function of $n$) is $n^{g/2} 2^{6g-7} n_g$ when $g$ is even and $\pi^{-1/2}n^{g/2} 2^{6g-7} m_g$ when $g$ is odd. The claim thus follows from the fact that $\gamma_k \sim (\pi k)^{-1/2}$ as $k\to \infty$. We then obtain the asymptotic expression of $L^+$ in ~\eqref{eq:L-fixed-g-large-n} by the asymptotic expressions of $v(g,n)$ and $u(g,n)$ in~\eqref{eq:v(g,n)-asy} and~\eqref{eq:u(g,n)-asy}.
\par 
Finally for $c_{\area}$, by~\eqref{eq:ekz} we have 
\begin{eqnarray*}
\frac{\pi^2}{3}c_{\area} & = & \frac{n+5-5g}{18} + L^{+}  \\
& \= & \frac{n+5-5g}{18} + \frac{1}{2g-3+n} \, \frac{r_g(n) + \gamma_{2g-3+n}s_g(n)}{p_g(n) + \gamma_{2g-3+n}q_g(n)}\,.  \\
\end{eqnarray*}
We thus define the polynomials  
$$ p^{*}_g(n) =  \frac{(n+5-5g)(2g-3+n)p_g(n)}{6} + 3r_g(n)\, $$
and
$$ q^{*}_g(n) =\frac{ (n+5-5g)q_g(n)}{6} + \frac{3s_g(n)}{2g-3+n}\,, $$
where $s_g(n)$ is divisible by $2g-3+n$ by definition. The above
expression for $c_{\area}$ can then be rewritten as
\begin{eqnarray*}
  c_{\area} & = &  \frac{1}{\pi^2} \,
  \cfrac{ \cfrac{ p^*_g(n)}{2g-3+n}  +   \gamma_{2g-3+n} q_g^*(n)}
        {p_g(n) + \gamma_{2g-3+n}q_g(n)}\,, 
\end{eqnarray*}
in accordance with~\eqref{eq:sv-fixed-g}. The claim about the degrees
of~$p^{*}_g(n)$
and~$q^{*}_g(n)$ follows from the degrees of $p_g(n)$ and $q_g(n)$.
\end{proof}

\subsection{$c_{\area}$ and $L^+$ in genus one}
\label{subsec:sv-L-g=1}

In this section we prove Corollary~\ref{cor:sv-L-g=1}. We have 
$$ L^+(1^n, -1^n) \= -2 \frac{ \int_{\BM_{1,n}}s_{n-1}(\BQ_{1,n})  \lambda_1}{\int_{\BM_{1,n}} s_{n}(\BQ_{1,n})}\,. $$
The denominator was computed previously in~\eqref{eq:segre-g=1}. For the numerator, we have 
\bas
\int_{\BM_{1,n}}s_{n-1}(\BQ_{1,n})\lambda_1 & \=  (-1)^{n-1} \int_{\BM_{1,n}} \left(\frac{\kappa_{(n-1)}\lambda_1}{(n-1)!} +  \frac{\kappa_{(n-2)}\lambda_1^2}{(n-1)!}\right) \\
& \=  (-1)^{n-1}\frac{(2n-3)!!}{24}\,, 
\eas
where $\kappa_{(n-1)}\lambda_1$ was computed previously in~\eqref{eq:segre-next-g=1}, and $\lambda_1^2 = 0$ as $\lambda_1$ is a pullback from $\BM_{1,1}$. It follows that 
\begin{eqnarray*} 
 L^+(1^n, -1^n)  \=  \frac{2}{1+ \frac{(2n-2)!!}{(2n-3)!!}}\,,
\end{eqnarray*}
thus proving~\eqref{eq:L-g=1}. Recall that $\frac{(2n)!!}{(2n-1)!!} \sim \sqrt{\pi n}$ for $n\to \infty$. Then~\eqref{eq:L-g=1} also implies a decay of $L^+$ in the order of $\frac{1}{\sqrt{n}}$ conjectured by Fougeron. 
\par
Finally by~\eqref{eq:ekz} we conclude that 
\bas 
 \frac{\pi^2}{3} c_{\area}(1^n, -1^n) & \=   L^+(1^n, -1^n) + \frac{n}{18} \\
 & \=    \frac{2}{1+ \frac{(2n-2)!!}{(2n-3)!!}}  + \frac{n}{18}\,,
\eas
thus proving~\eqref{eq:sv-g=1}. In particular, it confirms~\cite[Table 8]{TRMVV} for the case $g=1$ with $p_1^{*}(n) = \frac{1}{36}(n^2-n)$ and 
$q_1^{*}(n) = \frac{1}{36}n + 1$. 

\subsection{$c_{\area}$ and $L^+$ of affine invariant submanifolds}
\label{subsec:ais-sv-L}

Recall that a stratum of quadratic differentials can be lifted to the corresponding stratum of Abelian differentials via the canonical double cover, such that the image becomes an affine invariant submanifold. In general, let $\cNN$ be an arbitrary affine invariant submanifold in a stratum of Abelian differentials $\cHH(\mu)$. Suppose the tangent space of $\cNN$ projects onto a subspace $A$ of absolute periods, with kernel $R$ of relative periods.  
Denote $\dim_{\CC} A = a$ and $\dim_{\CC} R = r$, so that $\dim_{\CC} \cNN = a + r$. Without loss of generality, assume that a basis of $R$ is given by integration over $r$ paths joining the zeros $z_1,\ldots, z_r$ to a reference zero.   
\par 
Denote by $\PP\BNN$ the closure of $\PP\cNN$ in the IVC compactification of the projectivized stratum. Then 
$\dim_{\CC} \PP\BNN = (a-1) + r$.  Let $\xi$ be the first Chern class of the universal line bundle $\cOO(1)$ and $\delta$ the boundary divisor class. We make the following bold conjecture. 
\par 
\begin{conjecture}
\label{conj:ais-sv-L}
The area Siegel-Veech constant and sum of Lyapunov exponents of $\cNN$ can be obtained as the following intersection numbers: 
\ba
\label{eq:ais-sv}
 c_{\area}(\cNN) = -\frac{1}{4\pi^2} \frac{ \int_{\PP\BNN}\xi^{a-2} \psi_{1} \cdots \psi_{r} \delta }{ \int_{\PP\BNN}\xi^{a-1} \psi_{1} \cdots \psi_{r}}\,, 
\ea
\ba
\label{eq:ais-L}
 L(\cNN) = - \frac{ \int_{\PP\BNN}\xi^{a-2} \psi_{1} \cdots \psi_{r} \lambda_1}{ \int_{\PP\BNN}\xi^{a-1} \psi_{1} \cdots \psi_{r}}\,,
\ea
where $\psi_i$ is associated with the zero $z_i$ in the chosen basis of $R$. 
\end{conjecture}
\par 
We briefly explain the idea behind this conjecture. Since we work with the projectivized stratum, we can set one absolute period to be $1$, hence $\xi^{a-1}$ governs the 
absolute part of the volume form of $\cNN$, and $\psi_1, \ldots, \psi_r$ govern the relative part by varying the relevant zeros in $\cNN$. Therefore, the denominator in~\eqref{eq:ais-sv} and~\eqref{eq:ais-L} can be regarded as the volume of $\cNN$, up to a volume normalization factor. By now it has become clear that the boundary divisor class is responsible for $c_{\area}$ and the first Chern class of the Hodge bundle is responsible for $L$, thus explaining the structure of the conjectural formulas. In particular, using a different volume normalization factor should not matter, as it would cancel out between the numerator and denominator in each of the formulas. 
\par 
Moreover, there are evidences to support this conjecture from a number of known cases. If $\PP\cNN$ is a Teichm\"uller curve, then the conjectural formulas reduce to evaluating $\deg \delta / \deg \xi$ and $\deg \lambda_1 / \deg \xi$ up to the normalizing factors, and this case was well-understood after the works~\cite{kontsevich, ekz, chenrigid, cmNV}. If $\PP\cNN$ is a Hurwitz space of torus covers, then the above conjecture was established in~\cite[Section 4]{cmz}. If $\cNN$ is the entire stratum $\cHH(\mu)$, as said the above conjecture was verified in~\cite[Theorem 1.4]{CMSZ}. If $\cNN$ arises from a stratum of quadratic differentials via the canonical double cover, then Conjecture~\ref{conj:ais-sv-L} reduces to Conjecture~\ref{conj:sv-L}, which was proved in Theorem~\ref{thm:sv-L} for the case of the principal strata. Note that in this case $\zeta = 2\xi$ after lifting via the canonical double cover, which explains the difference of a factor $2$ in the two conjectures. 
\par 
We plan to treat the conjectures in this paper in future work.  

%% file: sec_appendix.tex
In \cite{TRMVV} we constructed a collection of generating series $(W_{g,n}^{{\rm I}})_{g \geq 0}^{n \geq 1}$ encoding some aspects of length statistics of multicurves, in which the Masur--Veech volumes of the principal strata of quadratic differentials appear as the lowest coefficients, and which satisfy the Eynard--Orantin topological recursion for a spectral curve $\mathcal{S}_{{\rm I}}$. This Appendix shows that a generating series $W_{g,n}^{{\rm II}}$ of intersection indices of the Segre class of Section~\ref{sec:Segre} with $\psi$-classes, in which the Masur--Veech volume is also the lowest coefficient, satisfy the same topological recursion for a (very different) spectral curve $\mathcal{S}_{{\rm II}}$. The two generating series have different meanings and are not a priori related. Only their lowest coefficients agree. We first review in Section~\ref{TRsummary} the definitions of the topological recursion, which originate in\cite{EORev}, in a simplified fashion which is sufficient for our needs. The main result of this Appendix is exposed in Proposition~\ref{secontr}. We prove it in Section~\ref{S3app} as a direct consequence (after some algebraic manipulations) of general relations between topological recursion and intersection theory established in \cite{EInter} and reviewed in Section~\ref{S2app}.

\section{Topological recursion for Masur--Veech volumes}
\label{TRsummary}

\subsection{Definition}
\label{TRdef}
For us, a spectral curve will be a quadruple $\mathcal{S} = (\mathcal{C},x,y,\omega_{0,2})$ as follows. $\mathcal{C}$ is an open subset of $\mathbb{P}^1$, $x$ is a (perhaps multivalued) function on $\mathcal{C}$ such that ${\rm d}x$ is meromorphic with a single, simple zero at $a \in \mathcal{C} \setminus \{\infty\}$, $y$ is a holomorphic function on $\mathcal{C}$ such that ${\rm d} y(a) \neq 0$, and $\omega_{0,2}$ is a meromorphic bidifferential whose only singularity on $\mathcal{C}^2$ is a double pole with biresidue $1$ on the diagonal.

\medskip

We define $\sigma$ to be the holomorphic involution defined in a neighbourhood $\mathcal{U} \subseteq \mathcal{C}$ of $a$ such that $\sigma(a) = a$, $x \circ \sigma = x$ and $\sigma \neq {\rm id}$. We introduce the recursion kernel
\[
	\mathfrak{K}(z_0,z) \,\coloneqq\, \frac{1}{2}\,\frac{\int_{\sigma(z)}^{z} \omega_{0,2}(\cdot,z_0)}{(y(z) - y(\sigma(z))){\rm d}x(z)}\,,
\]
which is a $1$-form in the variable $z_0 \in \mathcal{C}$ and a $(-1)$-form in the variable $z \in \mathcal{U}$. This allows the definition of multidifferentials $\omega_{g,n}$ on $\mathcal{C}^n$, indexed by $g \geq 0$ and $n > 0$ with $2g - 2 + n \geq 0$, by the following induction on $2g - 2 + n$:
\begin{flalign}
	\omega_{g,n}(z_1,\ldots,z_n)
	& \=
	\mathop{{\rm Res}}_{z = a} \; \mathfrak{K}(z_1,z) \bigg( \omega_{g - 1,n + 1}(z,\sigma(z),z_2,\ldots,z_n)\label{residuef} \\
	& \qquad\qquad\qquad
	+ \sum_{\substack{h + h' = g \\ J \sqcup J' = \{z_2,\ldots,z_n\}}} \omega_{h,1 + |J|}(z,J)\omega_{h',1 + |J'|}(\sigma(z),J') \bigg)\,, \nonumber
\end{flalign}
with the convention that $\omega_{0,1} = 0$. To be precise, $\omega_{g,n} \in H^0(\mathcal{C}^n,K_{\mathcal{C}}(*a)^{\boxtimes n})$ where $*a$ means allowing poles of arbitrary order at $a$, and although it is not apparent in their definition, $\omega_{g,n}$ are invariant under permutation of their $n$ variables. For $n = 0$ and $g \geq 2$, we also define the numbers
\begin{equation}
\label{Fg0}
	\omega_{g,0} \= \frac{1}{2 - 2g} \mathop{\rm Res}_{z = a} \bigg(\int_{a}^{z} y \, {\rm d}x\bigg)\omega_{g,1}(z)\,.
\end{equation}
We call $\omega_{g,n}$ the TR amplitudes.
\par
\medskip
The $\omega_{g,n}$ for $2g - 2 + n > 0$ can be decomposed
\begin{equation}
\label{omegagnbasis}
	\omega_{g,n}(z_1,\ldots,z_n)
	\=
	\sum_{\substack{k_1,\ldots,k_n \geq 0 \\ k_1 + \cdots + k_n \leq 3g - 3 + n}}
		F_{g,n}[k_1,\ldots,k_n] \, \prod_{i = 1}^n \xi_{k}(z_i)
\end{equation}
on the basis of $1$-forms $(\xi_k)_{k \geq 0}$ defined by
\[
	\xi_{0}(z_0) \,\coloneqq\, \mathop{\rm Res}_{z = a}\bigg(
		\frac{\omega_{0,2}(z_0,z)}{\sqrt{2(x(z) - x(a))}}
	\bigg)\,,
	\qquad
	\xi_{k} \,\coloneqq\, -{\rm d}\bigg(\frac{\xi_{k - 1}}{{\rm d}x}\bigg)\,.
\]
For $n = 0$ we also use the notation $F_{g,0} = \omega_{g,0}$ for uniformity.

\subsection{Applications to Masur--Veech volumes}

Here is the first topological recursion announced in the introduction of the Appendix.

\begin{Thm}[\cite{TRMVV}]
	Let $\omega_{g,n}^{{\rm I}}$ be the TR amplitudes for the spectral curve $\mathcal{S}_{{\rm I}}$ where $\mathcal{C}$ is a small neighbourhood of $0$ in $\mathbb{C}$, $x(z) = z^2/2$, $y(z) = -z$ and
	\[
		\omega_{0,2}^{{\rm I}}(z_1,z_2) \= \frac{{\rm d}z_1{\rm d}z_2}{2}\bigg(\frac{1}{(z_1 - z_2)^2} + \frac{\pi^2}{\sin^2\pi(z_1 - z_2)}\bigg)\,.
	\]
	For $2g - 2 + n > 0$ we have
        \bas \ 
		& \phantom{\=}\,\, {\rm vol}(\mathcal{Q}_{g,4g - 4 + 2n}(1^{4g - 4 + n},-1^{n})) \\ 
		& \=  \frac{2^{4g - 2 + n}(4g - 4 + n)!}{(6g - 7 + 2n)!}\,F_{g,n}^{{\rm I}}[0,\ldots,0] \\
		& \= \frac{2^{4g - 2 + n}(4g - 4 + n)!}{(6g - 7 + 2n)!} \mathop{\rm Res}_{z_1 = 0} \cdots \mathop{\rm Res}_{z_n = 0} \omega_{g,n}^{{\rm I}}(z_1,\ldots,z_n) \prod_{i = 1}^n z_i\,,
	\eas
where the third line is only valid for $n > 0$. For $n = 0$ and $g \geq 2$ we
have
$$
{\rm vol}(\mathcal{Q}_{g,4g - 4}(1^{4g - 4},-1^0)) \= \frac{3\cdot 2^{4g - 2}(4g - 4)!}{(6g - 6)!}\,F_{g,1}^{{\rm I}}[1]\,.
$$
\end{Thm}
\par
In this Appendix we show a second topological recursion.
\par
\begin{Prop}
\label{secontr}
	Let $\omega_{g,n}^{{\rm II}}$ be the TR amplitudes for the spectral curve $\mathcal{S}_{{\rm II}}$ defined by
	\[
		\mathcal{C} = \mathbb{P}^1\,,\qquad x(z) = -z - \ln z\,,\qquad y(z) = z^2\,,
                \qquad \omega_{0,2}(z_1,z_2) = \frac{{\rm d}z_1{\rm d}z_2}{(z_1 - z_2)^2}\,.
	\]
	For $2g - 2 + n > 0$ and $k_1,\ldots,k_n \geq 0$ we have
	\begin{equation}
	\label{firstunf}
		F_{g,n}^{{\rm II}}[k_1,\ldots,k_n] = 2^{2 - 2g - n} \int_{\overline{\mathcal{M}}_{g,n}} s(\overline{\mathcal{Q}}_{g,n}) \prod_{i = 1}^n \psi_i^{k_i}\,.
	\end{equation}
	In particular, in view of Proposition~\ref{prop:volume-segre},
        we have 
	\begin{flalign}
		& \phantom{\=}\,\, {\rm vol}(\mathcal{Q}_{g,4g - 4 + 2n}(1^{4g - 4 + n},-1^{n})) \nonumber \\
		& = \frac{2^{4g - 1 + n}({\rm i}\pi)^{6g - 6 + 2n}(4g - 4 + n)!}{(6g - 7 + 2n)!}\,F_{g,n}^{{\rm II}}[0,\ldots,0] 	\label{secondunf} \\
		& = \frac{2^{4g - 1 + n}({\rm i}\pi)^{6g - 6 + 2n}(4g - 4 + n)!}{(6g- 7 + 2n)!} \mathop{\rm Res}_{z_1 = -1} \cdots \mathop{\rm Res}_{z_n = -1} \omega_{g,n}^{{\rm II}}(z_1,\ldots,z_n) \prod_{i = 1}^n (z_i + 1)\,,
                \nonumber
	\end{flalign}
	where the third line is only valid for $n > 0$. For $n = 0$ and $g \geq 2$ we have
	\begin{equation}
	\label{thenocase}
		{\rm vol}(\mathcal{Q}_{g,4g - 4}(1^{4g - 4},-1^{0})) \= \frac{3\cdot 2^{4g}({\rm i}\pi)^{6g - 6}(4g - 4)!}{(6g - 6)!}\big(F_{g,1}^{{\rm II}}[1] + F_{g,1}^{{\rm II}}[2]\big)\,.
	\end{equation}
\end{Prop}

\section{Topological recursion and intersection theory}
\label{S2app}

Let $\mathcal{S}$ be a spectral curve as in Section~\ref{TRdef} and $\omega_{g,n}$ the corresponding TR amplitudes. The coefficients in \eqref{omegagnbasis} can then be interpreted in terms of intersection theory on $\overline{\mathcal{M}}_{g,n}$. To state the formula, we introduce two formal power series:
\begin{align}
	\label{Teqn}
	T(u) & \= \frac{e^{x(a)u^{-1}}}{\sqrt{2\pi u}} \int_{\gamma} e^{-x(z)u^{-1}}{\rm d}y(z) \= \exp\bigg(-\sum_{d \geq 0} t_{d}\,u^d\bigg)\,, \\
	\label{Reqn}
	R(u) & \= \frac{e^{x(a)u^{-1}}}{\sqrt{2\pi u^{-1}}} \int_{\gamma} e^{-x(z)u^{-1}} \xi_0(z) \= \exp\bigg(\sum_{d \geq 1} r_{d}\,u^d\bigg) = \exp\bigl( r(u) \bigr)\,. 
\end{align}
Here, $\gamma$ is the steepest descent contour for the function $x/u$ on $\mathcal{C}$, going around $a$ in the positive direction. $T(u)$ and $R(u)$ are the asymptotic series to the right-hand side when $u \rightarrow +\infty$, and these definitions only depend on the germ of $\gamma$ near $a$ (see below).

\medskip

We define a class $\Omega_{g,n} \in H^{\bullet}(\overline{\mathcal{M}}_{g,n})$ by the formula
\begin{flalign}
	\Omega_{g,n}
	\=
	\exp\bigg(
	&
	\sum_{d \geq 0} t_{d} \, \kappa_{d} + \sum_{i = 1}^n r(\psi_i) \label{theclasss} \\
	&
	+
	\frac{1}{2} i_{{\rm irr }\ast}
		\bigg(\frac{r(\psi_{1}) + r(\psi_{2})}{\psi_{1} + \psi_{2}}\bigg)
	+
	\frac{1}{2} \sum_{h=0}^g \sum_{S \subseteq [[1,n]]} i_{{h,S}\ast}
		\bigg(\frac{r(\psi_{1}) + r(\psi_{2})}{\psi_{1} + \psi_{2}}\bigg) 
	\bigg)\,, \nonumber
\end{flalign}
where $i_{{\rm irr}}$ and $i_{h,S}$ are the maps introduced in Section~\ref{sec:Segre}.
\begin{Thm}
\label{intrep}
	For $2g - 2 + n > 0$, we have the equality in $\mathbb{C}[[\mu_1^{-1},\ldots,\mu_n^{-1}]]$
	\begin{equation}
	\label{ints}
		\int_{\gamma^n} \omega_{g,n}(z_1,\ldots,z_n) \, \prod_{i = 1}^n \frac{e^{- \mu_i (x(z_i) - x(a))}{\rm d}z_i}{\sqrt{2\pi\mu_i}\,R(1/\mu_i)} \= \int_{\overline{\mathcal{M}}_{g,n}} \Omega_{g,n}\,\, \prod_{i = 1}^n \frac{1}{1 + \mu_i\psi_i}\,.
	\end{equation}
	Equivalently, in the decomposition \eqref{omegagnbasis} we have for $k_1,\ldots,k_n \geq 0$
	\bes
		F_{g,n}[k_1,\ldots,k_n] \= \int_{\overline{\mathcal{M}}_{g,n}} \Omega_{g,n}\,\, \prod_{i = 1}^n \psi_i^{k_i}\,.
	\ees
\end{Thm}

\begin{proof}
	We explain how to derive the particular form we give to the result \eqref{ints} from \cite[Theorem 3.1]{EInter}.

	\smallskip

	Firstly, the Laplace variable $u$ in \cite{EInter} is the variable $u^{-1}$ for us. We chose this convention, as we found more convenient to work with formal power series in~$u$ (instead of~$u^{-1}$).

	\smallskip

	Secondly, as it is clear from the proof in \cite{EInter}, the contribution of boundary divisors in \cite[Formula~(3.11)]{EInter} should not be understood as the genuine exponential of a class, but rather as a sum over stable graphs, where the weight of the edges is $\check{B}(1/\psi_{1},1/\psi_{2})$, which we here denote $\mathfrak{E}(\psi_{1},\psi_{2})$. $\mathfrak{E}(u_1,u_2)$ is a formal power series in $u_1$ and $u_2$ which can be computed form the data of $x$ and $\omega_{0,2}$. Since ${\rm d}x$ is meromorphic on the compact curve $\mathbb{P}^1$, we can use \cite[Appendix~B]{EInter} which justifies (with the preceding conventions) that
	\begin{equation} 
	\label{relE}
		\mathfrak{E}(u,v) = \frac{1 - R(u)R(v)}{u + v}\,.
	\end{equation}
	Our $R(u)$ corresponds to $\pm f_{a,0}(u^{-1})$ in \cite{EInter}. The sign depends on the choice of square root, that should be made so that $R(u) = 1 + O(u)$, but  it does not affect \eqref{relE} since $R$ appears by pairs. We also warn the reader familiar with cohomological field theories that $R$ rather corresponds to the inverse of the $R$-matrix in Givental formalism.

	\smallskip

	Thirdly, the sum over stable graphs can be converted into intersections in $\overline{\mathcal{M}}_{g,n}$ of the exponential of a boundary class. Namely \cite[Formula (3.11)]{EInter} is correctly interpreted as involving the exponential of a boundary class if we replace the contribution of boundary divisors by the pushforward of $E(\psi_{1},\psi_{2})$, where the new generating series is
	\[
		E(u_1,u_2) \= -\frac{\ln\big(1 - (u_1 + u_2)\mathfrak{E}(u_1,u_2)\big)}{u_1 + u_2}\,,
	\]
	or equivalently
	\[
		\mathfrak{E}(u_1,u_2) \= \frac{1 - e^{-(u_1 + u_2)E(u_1,u_2)}}{u_1 + u_2}\,.
	\] 
	This relation comes from taking into account self-intersections of divisors, see e.g.\ \cite[Lemma 3.10]{TRMVV}.

	\smallskip

	Fourthly, the relation \eqref{relE} leads to a simplification of the contribution of $\psi$-classes in \cite[Formula~(3.10)]{EInter}, namely
	\[
		\frac{2\sqrt{\pi}\,e^{-\mu_i x(a)}}{\sqrt{\mu_i}}\bigg(\frac{\mu_i}{1 + \mu_i\psi_i} - \mathfrak{E}(1/\mu_i,\psi_i)\bigg) \= 2\sqrt{\pi\mu_i}\,e^{-\mu_i x(a)} \frac{R(1/\mu_i)R(\psi_i)}{1 + \mu_i\psi_i}\,.
	\] 
	If we factor out $\prod_{i = 1}^n \sqrt{2\pi \mu_i}\,R(1/\mu_i)$ to put it in the left-hand side of \cite[Formula~(3.10)]{EInter}, this leaves in the right-hand side a power $2^{n/2}$, which combines with the overall power $2^{3g - 3 + n}$ to give $2^{(3/2)(2g - 2 + n)} = 2^{(3/2)\kappa_0}$. Therefore, we can change the definition of the coefficient of $\kappa_0$ in \cite{EInter}; it was there denoted $\hat{t}_0$ and it is related to our $t_0$ by $\hat{t}_0 + (3/2)\ln 2 = t_0$. Note that there seems to be a misprint in \cite{EInter}, where Formulas~(3.12) and (4.15) should have the prefactor of $2$ in the numerator rather than in the denominator. Making this correction led us to the definition of $T(u)$ with a prefactor of $2^{-1/2}$.
\end{proof}

For certain $(x,y)$, $T(u)$ and $R(u)$ can be identified with well-known special functions and $t_k$ and $r_k$ can be computed explicitly. A term-by-term computation is always possible, for instance as follows. Let $\zeta(z) = \sqrt{2(x(z) - x(a))}$ be the local coordinate near $a$ (for the standard determination of the square root), which has the property that $\zeta(\sigma(z)) = -z$, and compute the expansion near $z = a$:
\[
	y(z) \= \sum_{k \geq 0} y_{k}\,\zeta(z)^{k}\,,\qquad \xi_{0}(z) \= \frac{{\rm d}\zeta(z)}{\zeta(z)^2} + \sum_{k \geq 0} \xi_{0,k} \, \zeta(z)^{k} {\rm d}\zeta(z)\,.
\]
Then, to obtain the asymptotic expansion of the integrals up to $O(u^{\infty})$, we can take $u > 0$ and replace $\gamma$ in the $z$-plane with a contour $\delta$ in the $\zeta$-plane that goes from $ +\infty - {\rm i}0$ to $-{\rm i}0$, then follows the half-circle leaving $0$ to its right until $ {\rm i}0$, from where it goes to $+\infty + {\rm i}0$. We have
\[
	\int_{\delta} e^{-\zeta^2/2u}\,\zeta^{2k} \frac{{\rm d}\zeta}{\sqrt{2\pi u}} = -(2k - 1)!! u^{k}\,,
\]
where the global minus sign comes from the orientation of the contour. This formula remains valid for any $k \geq -1$, with the convention that $(-1)!! = 1$ and $(-3)!! = -1$. Thus
\[
	T(u) \= -\sum_{d \geq 0} (2d + 1)!!\,y_{2d + 1}\,u^{d}\,,\qquad R(u) \= 1 + \sum_{d \geq 0} (2d - 1)!!\,\xi_{0,2d}\,u^{d + 1}\,.
\]

\section{Study of a family of spectral curves}

\subsection{Definition and basic properties}
\label{S3app}
Let $(a,b) \in \mathbb{C}^{\ast} \times \mathbb{Z}^{\ast}$ and consider the spectral curve $\mathcal{S}[a,b]$ defined by
\begin{equation}
\label{spx}
	x(z) \= -z + a\ln z\,,
	\qquad
	y(z) \= z^{b}\,,
	\qquad
	\omega_{0,2}(z_1,z_2) \= \frac{{\rm d}z_1{\rm d}z_2}{(z_1 - z_2)^2}\,.
\end{equation}
Notice that ${\rm d}x$ has a unique, simple zero at $z = a$. To complete the definition, we choose $\mathcal{C}$ to be a small neighbourhood of $a$ in $\mathbb{C}$. The determination of the logarithm is chosen arbitrarily and will not affect our discussion. For the record, we compute
\[
	x(a) \= a\big(\ln(a) - 1\big)\,,\qquad x''(a) \= -\frac{1}{a}\,.
\]
We will use instead of $z$ the coordinate $t = z - a$, so that
\[
	x(z) \= -(t + a) + a\ln(t + a)\,,\qquad y(z) \= (t + a)^{b}\,.
\]
The involution such that $x(a + \sigma(t)) = x(a + t)$ is given by $\sigma(t) = -a^{-1}\hat{\sigma}(-a^{-1}t)$ where $\hat{\sigma}$ is the unique solution to
\[
	t - \hat{\sigma}(t) \= \ln\bigg(\frac{1 - \hat{\sigma}(t)}{1 - t}\bigg)\,,\qquad
        \hat{\sigma}(t) \= -t + O(t^2)\,.
\]
It does not have a simple expression, but can be generated to high order on the computer
\[
	\hat{\sigma}(t) = -\bigg(t + \frac{2t^2}{3} + \frac{4t^3}{9} + \frac{44t^4}{135} + \frac{104t^5}{405} + \frac{40t^6}{189} + \frac{7648t^7}{42525} + \frac{2848t^8}{18225} + \frac{31712t^{9}}{229635} + O(t^{10})\bigg)\,.
\]

\begin{Lemma}
	For $2g - 2 + n > 0$ and $k_1,\ldots,k_n \geq 0$ we have
	\[
		F_{g,n}[k_1,\ldots,k_n] \in ((-a)^{b - 1/2}b)^{2 - 2g - n}\cdot \mathbb{Q}[a^{-1},b]\,.
	\]
\end{Lemma}

\begin{proof}
	As we need to stress the dependence in $a$ and $b$ in this proof, we momentarily denote $\omega_{g,n}^{[a,b]}$ the TR amplitudes associated with $\mathcal{S}[a,b]$. We claim that for $2g - 2 + n > 0$ we have
	\begin{equation}
	\label{depa}
		\omega_{g,n}^{[a,b]}(z_1,\ldots,z_n)
		\=
		(-a)^{(b + 1)(2 - 2g - n)}\,\omega_{g,n}^{[-1,b]}(-z_1/a,\ldots,-z_n/a)\,.
	\end{equation}
	This is justified by noticing that $x(z) = -a\tilde{x}(-z/a)$ and $y(z) = (-a)^{b}(-z/a)^b$ with $\tilde{x}(\tilde{z}) = c + x(\tilde{z})|_{a = -1}$ for some constant $c$, $\tilde{y}(\tilde{z}) = y(\tilde{z})|_{a = -1}$ and $\omega_{0,2}(z_1,z_2) = \omega_{0,2}(-z_1/a,-z_2/a)$. Since $(x,y)$ are involved in \eqref{residuef} only via the $1$-form $y{\rm d}x$ in the denominator of the recursion kernel and $\omega_{g,n}$ is reached by $2g - 2 + n$ steps of the recursion, we deduce \eqref{depa} for $n > 0$. Inserting this result for $\omega_{g,1}$ in \eqref{Fg0}, we see that \eqref{depa} also holds for $n = 0$.

	We denote $\xi_k^{[a]}$ the basis of $1$-forms, since it only depends on $a$. One easily checks by induction on $k$ that
	\[
		\xi^{[a]}_{k}(z) = -(-a)^{-(k + 1/2)}\,\xi^{[-1]}_{k}(-z/a)\,.
	\]
	We deduce that
	\[
		F_{g,n}^{[a,b]}[k_1,\ldots,k_n]
		\=
		(-1)^{n} (-a)^{(b + 1)(2 - 2g - n) + \sum_{i = 1}^n (k_i + 1/2)}\,F_{g,n}^{[-1,b]}[k_1,\ldots,k_n]\,.
	\]
	Since in \eqref{omegagnbasis} we have $\sum_i k_i \leq 3g - 3 + n$, we obtain that, for any value of $b \in \mathbb{Z}^{\ast}$,
	\begin{equation}
	\label{depa2}
		F_{g,n}^{[a,b]}[k_1,\ldots,k_n]
		\in
		(-a)^{(b - 1/2)(2 - 2g - n)} \cdot \mathbb{Q}_{3g - 3 + n}[a^{-1}]\,.
	\end{equation}
	We now study the dependence in $b$. Taking into account the dependence in $a$ of the involution, we observe that the $t \rightarrow 0$ expansion of the recursion kernel $\mathfrak{K}(a + t_1,a + t)$ belongs to $t^{-1}(ba^{b - 2})^{-1} \cdot \mathbb{Q}[t_1^{-1},a^{-1},b][[t]]\cdot \frac{{\rm d}t_1}{{\rm d}t}$. This implies by induction that for $2g - 2 + n > 0$ and $n > 0$
	\begin{equation}
	\label{bdependence}
		\omega_{g,n}(t_1,\ldots,t_n) \in (ba^{b - 2})^{2 - 2g - n} \cdot \mathbb{Q}[a^{-1},b][t_1^{-1},\ldots,t_n^{-1}] \prod_{i = 1}^n \frac{{\rm d}t_i}{t_i}\,.
	\end{equation}
	Combining with \eqref{Fg0} then extends the validity of \eqref{bdependence} to $n = 0$, which together with \eqref{depa2} proves the claim.
\end{proof}

\subsection{Intersection theory}

We recall the Hankel representation for $v \in \mathbb{C}$
\[
	\frac{1}{\Gamma(v)} = \frac{1}{2{\rm i}\pi} \int_{c} e^{t}\,t^{-v} \, {\rm d}t\,,
\]
where the contour $c$ (as given by Theorem~\ref{intrep}) goes in the $v$-plane from $ -\infty - {\rm i}0$ to $-{\rm i}0$, then follows the half-circle leaving $0$ to its left until $ {\rm i}0$, from where it goes to $-\infty + {\rm i}0$. We also recall that for fixed $\beta \in \mathbb{C}$, we have the asymptotic expansion when $v \rightarrow \infty$ such that $|{\rm arg}\,v| < \pi - \epsilon$ for some fixed $\epsilon > 0$
\[
	\Gamma(v + \beta)
	\=
	\exp\bigg(v \ln v - v + \Big(\beta - \frac{1}{2}\Big)\ln v + \frac{\ln(2\pi)}{2} + \sum_{d \geq 1} \frac{(-1)^{d + 1}B_{d + 1}(\beta)}{d(d + 1)}\,v^{-d}\bigg)\,,
\]
where $B_d(\beta)$ are the Bernoulli polynomials and $B_d(0) = B_d$ are the Bernoulli numbers (see Section~\ref{sec:Segre}). They vanish if $d$ is odd and greater than $2$.

\medskip

We can compute $T(u)$ and $R(u)$ for the spectral curve \eqref{spx}, first setting $u > 0$ and $a \notin \mathbb{R}_{-}$.
\begin{equation*}
\begin{split}
T(u) & \= \frac{b}{\sqrt{2\pi u}}\,e^{au^{-1}(\ln(a) - 1)} \int_{\gamma} e^{zu^{-1}}z^{-au^{-1}} z^{b - 1} \,{\rm d}z \\
	& \= \frac{b}{\sqrt{2\pi u}} e^{au^{-1}(\ln(au^{-1}) - 1) - b\ln(u^{-1})} \int_{c} e^{t}t^{-au^{-1} + b- 1} \, {\rm d}t \\  
	& \= \frac{{\rm i}b\,e^{au^{-1}\ln(au^{-1}) - au^{-1} + (1/2 - b)\ln(u^{-1}) + \frac{1}{2}\ln(2\pi)}}{\Gamma(au^{-1} + 1 - b)}\,,
\end{split} 
\end{equation*}
where we used the Hankel representation. Using the asymptotic expansion for $\ln(\Gamma)$, we obtain when $u \rightarrow 0$
\[
	T(u) \= {\rm i}ba^{b - 1/2} \exp\bigg(\sum_{d \geq 1} \frac{(-1)^{d}\,B_{d + 1}(1 - b)}{d(d + 1) a^d}\,u^{d}\bigg)\,,
\]
that is
\[
	e^{t_0} \= -{\rm i}b^{-1}a^{1/2 - b}\,,
	\qquad
	t_d \= (-1)^{d + 1} \frac{B_{d + 1}(1 - b)}{d(d + 1)a^{d}}
	\qquad(\text{for }d \geq 1)\,.
\]
We compute using integration by parts
\[
\begin{aligned}
	R(u)
	& \= \frac{e^{au^{-1}(\ln(a) - 1)}}{\sqrt{-2\pi (au)^{-1}}} \int_{\gamma} \frac{e^{zu^{-1}}\,z^{-au^{-1}}}{(z - a)^2} \, {\rm d}z \\
	& \= \frac{e^{au^{-1}(\ln(a) - 1)}}{\sqrt{-2\pi (au^{-1})^{-1}}} \int_{\gamma} e^{zu^{-1}}z^{-au^{-1} - 1} \, {\rm d}z\,,
\end{aligned}
\]
where we notice the cancellation of poles between the integrated factor $(z - a)^{-1}$ and the factor coming from the derivative, the only effect being an extra factor of $(zu)^{-1}$ turning $u^{1/2}$ into $u^{-1/2}$ and $z^{-au^{-1}}$ into $z^{-au^{-1} - 1}$. We then get
\[
	R(u) \= \frac{\sqrt{2\pi au^{-1}}e^{au^{-1}(\ln(au^{-1}) - 1)}}{\Gamma(au^{-1} + 1)}\,.
\]
It admits the asymptotic expansion when $u \rightarrow 0$
\[
	R(u) \= \exp\bigg(-\sum_{d \geq 1} \frac{B_{d + 1}}{d(d + 1)a^{d}}\,u^{d}\bigg)\,.
\]
Therefore, we are led to define the class
\bas
\Omega_{g,n}[a,b]
	\,\coloneqq\,
		\exp\bigg\{
		\sum_{d \geq 1} \bigg( &
			\frac{(-1)^{d+1}B_{d + 1}(1 - b)}{d(d + 1)a^{d}}\,\kappa_d
			-
			\sum_{i = 1}^n \frac{B_{d + 1}}{d(d + 1)a^d}\,\psi_i^{d}
			\\
			&
			+
			\frac{1}{2} \frac{B_{d + 1}}{d(d + 1)a^{d}}
			\sum_{i+j = d-1} i_{{\rm irr}\ast}
				(-\psi_1)^i \psi_2^j \\
			&
			+
			\frac{1}{2} \frac{B_{d + 1}}{d(d + 1)a^{d}}
			\sum_{h=0}^g \sum_{S \subseteq [[1,n]]} \sum_{i+j = d-1}
				i_{{h,S}\ast} (-\psi_1)^i \psi_2^j
		\bigg)\bigg\}\,,
\eas
where compared to \eqref{theclasss} we have excluded the $\kappa_0$ term. Due to the vanishing of odd Bernoulli numbers and the symmetry $(-1)^{d}B_{d}(1 - b) = B_{d}(b)$, we have $\Omega_{g,n}[-a,b] = \Omega_{g,n}^{-1}[a,1 - b]$. Notice that $\Omega_{g,n}[a,b]$ is a polynomial in $a^{-1}$ (for each fixed $g,n$), therefore makes sense even if~$a$ is on the negative real axis.
 
\begin{Cor}
\label{Co02}
	For any $(a,b) \in \mathbb{C}^{\ast} \times \ZZ^{\ast}$, the TR amplitudes of the spectral curve \eqref{spx} are decomposed as in \eqref{omegagnbasis} with
	\begin{equation}
	\label{fgndfs}
		F_{g,n}[k_1,\ldots,k_n] \= ({\rm i}ba^{b - 1/2})^{2 - 2g - n}\,\int_{\overline{\mathcal{M}}_{g,n}} \Omega_{g,n}[a,b] \prod_{i = 1}^n \psi_i^{k_i}\,.
	\end{equation}
	In particular,
	\bes
		F_{g,n}[0,\ldots,0] \= ({\rm i}ba^{b - 1/2})^{2 - 2g - n} \int_{\overline{\mathcal{M}}_{g,n}} \Omega_{g,n}[a,b]\,.
	\ees
\end{Cor}

\begin{proof}
	For $a \notin \mathbb{R}_{-}$, we just apply Theorem~\ref{intrep}. We know a priori that the left-hand side of \eqref{fgndfs} divided by $(ba^{b - 1/2})^{2 - 2g - n}$ is a polynomial in $a^{-1}$ and $b$, and the class $\Omega_{g,n}[a,b]$ depends polynomially on $a^{-1}$ and $b$. Therefore, the equality holds for all $(a,b) \in \mathbb{C}^{\ast} \times \mathbb{Z}^{\ast}$.
\end{proof}

\begin{Rem}
	The Chiodo class ${\rm ch}(R^{\bullet}\pi_*\mathcal{S})$ of \cite{Chiodo} for the values $(r,s) = (1,b)$ coincides with our class $\Omega_{g,n}[1,b]$. For $b = 1$ this is simply the Chern character of the Hodge bundle: the first proof of Corollary~\ref{Co02} in that case comes from combining \cite{ELSV} (ELSV formula for Hurwitz numbers) and topological recursion for Hurwitz numbers \cite{EMS}; computations similar to ours appear in \cite{SSZ}.
\end{Rem}

\subsection{Specialization to $(a,b) = (-1,2)$}

For these values, the class we have constructed is precisely the Segre class appearing
in Lemma~\ref{lm:chern}.
\bas \
	\Omega_{g,n}[-1,2]
	& \=
		\exp\bigg(
		\sum_{d \geq 1} \bigg(
			\frac{(-1)^d B_{d + 1}(2)}{d(d + 1)}\,\kappa_d
			-
			\sum_{i = 1}^n \frac{(-1)^d B_{d + 1}}{d(d + 1)}\,\psi_i^{d}
			\\
		& \qquad\qquad\qquad\quad
			+
			\frac{1}{2} \frac{(-1)^d B_{d + 1}}{d(d + 1)}
			\sum_{i+j = d-1} i_{{\rm irr }\ast}
				(-\psi_1)^i \psi_2^j \\
		& \qquad\qquad\qquad\quad
			+
			\frac{1}{2} \frac{(-1)^d B_{d + 1}}{d(d + 1)}
			\sum_{h=0}^g \sum_{S \subseteq [[1,n]]} \sum_{i+j = d-1}
				i_{{h,S}\ast} (-\psi_1)^i \psi_2^j
		\bigg) \\
	& \=
		s(\overline{\mathbb{E}}_{g,n})
		\exp\bigg( \sum_{d \geq 1} \frac{(-1)^d}{d} \kappa_d \bigg) \\
	& \= s(\overline{\mathcal{Q}}_{g,n})\,.
\eas
The application of Corollary~\ref{Co02} yields for $2g - 2 + n > 0$ and $k_1,\ldots,k_n \geq 0$
\begin{equation}
\label{fgnseg}
	F_{g,n}[0,\ldots,0] \= 2^{2 - 2g - n}\,\int_{\overline{\mathcal{M}}_{g,n}} s(\overline{\mathcal{Q}}_{g,n})\,.
\end{equation}
Recalling Proposition~\ref{prop:volume-segre}, we have proved
Theorem~\ref{secontr}, except for the alternative formula \eqref{thenocase}
for the $n = 0$ case. To obtain it, we go back to \eqref{Fg0} which gives
\begin{equation}
\label{ozerres}
	F_{g,0} \= \frac{1}{2g - 2} \mathop{\rm Res}_{z = -1} \frac{(2z - 1)(z + 1)^2}{6} \bigg(\sum_{k = 0}^{3g - 3} F_{g,1}[k]\,\xi_{k}(z)\bigg)\,.
\end{equation}
It is easy to prove by induction on $k \geq 0$ that
\[
	\xi_{k}(z) \= \sum_{l = k}^{2k} \frac{c_{k,l}\,{\rm d}z}{(z + 1)^{l + 2}}\,,\qquad c_{k,l} \in \mathbb{Z}\,.
\]
Therefore, only the terms $k \in \{1,2\}$ contribute to the residue in \eqref{ozerres} and we find
\[
	F_{g,0} \= \frac{F_{g,1}[1] + F_{g,1}[2]}{g - 1}\,,
\]
which can be rearranged into the desired Equation~\eqref{thenocase}.

For computations of the TR amplitudes, it is simpler to work with rational functions instead of rational differential forms. We therefore set
\bas
	W_{g,n}(t_1,\ldots,t_n)
	& \= \frac{\omega_{g,n}(t_1 + a,\ldots,t_n + a)}{{\rm d}t_1 \cdots {\rm d}t_n}\,, \\
	K(t_1,t) & \= \mathfrak{K}(t_1 + a,t + a) \sigma'(t + a)\,\frac{{\rm d}t}{{\rm d}t_1}\,, \\
	\Xi_k(t) & \= \frac{\xi_k(t)}{{\rm d}t}\,.
\eas
Then, $W_{g,n}(t_1,\ldots,t_n) \in \prod_{i = 1}^n t_i^{-1} \cdot \mathbb{C}[t_1^{-1},\ldots,t_n^{-1}]$. We focus on the case $(a,b) = (-1,2)$, in which case
\bes
	K(t_1,t) \= \frac{t\hat{\sigma}'(t)}{2(t + 1)(2 - t - \hat{\sigma}(t))}\,\frac{1}{(t_1 - t)(t_1 - \hat{\sigma}(t))}\,.
\ees
The recursion formula becomes
\ba \label{wngfun}
W_{g,n}(t_1,\ldots,t_n) 
& \= \mathop{\rm Res}_{t = 0} {\rm d}t\,K(t_1,t)\bigg(W_{g - 1,n + 1}(t,\hat{\sigma}(t),t_2,\ldots,t_n) \\
        & \phantom{\=} \quad + \sum_{\substack{J \sqcup J' = \{t_2,\ldots,t_n\} \\ h + h' = g}} \!\!\!\! W_{h,1 + |J|}(t,J)W_{h',1 + |J'|}(\hat{\sigma}(t),J')\bigg)\,.
\ea
The first elements on the basis in which we can read the $F_{g,n}$ are
\begin{equation*}
	\begin{aligned}
		\Xi_0(t) & \= \tfrac{1}{t^2} \\
		\Xi_1(t) & \= -\tfrac{2}{t^3} + \tfrac{3}{t^4} \\
		\Xi_2(t) & \= \tfrac{6}{t^4} - \tfrac{20}{t^5} + \tfrac{15}{t^6}
	\end{aligned}
	\qquad 
	\begin{aligned}
		\Xi_3(t) & \= -\tfrac{24}{t^5} + \tfrac{130}{t^6} - \tfrac{210}{t^7} + \tfrac{105}{t^8} \\
		\Xi_4(t) & \= \tfrac{120}{t^6} - \tfrac{924}{t^7} + \tfrac{2380}{t^8} - \tfrac{2520}{t^9} + \tfrac{945}{t^{10}} \\
		\Xi_5(t) & \= -\tfrac{720}{t^7} + \tfrac{7308}{t^8} - \tfrac{26432}{t^9} + \tfrac{44100}{t^{10}} - \tfrac{34650}{t^{11}} + \tfrac{10395}{t^{12}} 
	\end{aligned}
\end{equation*}
We need to expand the recursion kernel when $t \rightarrow 0$
\[
	K(t_1,t) = \sum_{j \geq -1} K_j(t_1)\,t^{j}
\]
and it is useful to decompose the coefficients on the $(\Xi_m(t_1))_{m \geq 0}$.
\begin{equation*}
	\begin{aligned}
		K_{-1} & \= \tfrac{\Xi_0}{4} \\
		K_{0} & \= \tfrac{\Xi_0}{12} \\
		K_1 & \= -\tfrac{\Xi_0}{12} + \tfrac{\Xi_1}{12} \\
		K_2 & \= -\tfrac{49\,\Xi_0}{540} + \tfrac{\Xi_1}{12}
	\end{aligned}
	\qquad\quad\qquad
	\begin{aligned}
		K_3 & \= -\tfrac{59\,\Xi_0}{1620} + \tfrac{17\,\Xi_1}{540} + \tfrac{\Xi_2}{60} \\
		K_4 & \= -\tfrac{\Xi_0}{2268} - \tfrac{\Xi_1}{324} + \tfrac{\Xi_2}{36} \\
		K_5 & \= \tfrac{1021\,\Xi_0}{170100} - \tfrac{11\,\Xi_1}{1260} + \tfrac{97\,\Xi_2}{3780} + \tfrac{\Xi_3}{420} \\
		K_6 & \= \tfrac{17\,\Xi_0}{72900} - \tfrac{59\,\Xi_1}{24300} + \tfrac{149\,\Xi_2}{8100} + \tfrac{\Xi_3}{180} \\
	\end{aligned}
\end{equation*}
Applying \eqref{wngfun} and rearranging the result as a multilinear combination of $\Xi_{k_i}(t_i)$ we arrive to Table~\ref{tablak}.
\begin{table}[h!]
\centering
\begin{tabular}{cc}
\begin{tabular}[t]{|c|c|c|}
\hline
$(g,n)$ {\rule{0pt}{3.2ex}}{\rule[-1.8ex]{0pt}{0pt}}& $(k_1,\ldots,k_n)$ & $F_{g,n}[\mathbf{k}]$ \\
\hline
\hline
$(0,3)$ {\rule{0pt}{3.2ex}}{\rule[-1.8ex]{0pt}{0pt}}& $(0,0,0)$ & $\frac{1}{2}$ \\
\hline
\multirow{2}{*}{$(0,4)$} {\rule{0pt}{3.2ex}}{\rule[-1.8ex]{0pt}{0pt}}& $(0,0,0,0)$ & $\mathbf{-\frac{1}{4}}$ \\
{\rule{0pt}{3.2ex}}{\rule[-1.8ex]{0pt}{0pt}}& $(1,0,0,0)$ & $\frac{1}{4}$ \\
\hline 
\multirow{4}{*}{$(0,5)$}  {\rule{0pt}{3.2ex}}{\rule[-1.8ex]{0pt}{0pt}}& $(0,0,0,0,0)$& $\mathbf{\frac{3}{8}}$ \\
{\rule{0pt}{3.2ex}}{\rule[-1.8ex]{0pt}{0pt}}&  $(1,0,0,0,0)$ & $-\frac{3}{8}$ \\
{\rule{0pt}{3.2ex}}{\rule[-1.8ex]{0pt}{0pt}}& $(2,0,0,0,0)$ & $\frac{1}{8}$ \\
{\rule{0pt}{3.2ex}}{\rule[-1.8ex]{0pt}{0pt}}& $(1,1,0,0,0)$ & $\frac{1}{4}$ \\
\hline
\hline
\multirow{5}{*}{$(2,1)$}  {\rule{0pt}{3.2ex}}{\rule[-1.8ex]{0pt}{0pt}}& $(0)$ & $\mathbf{\frac{29}{5120}}$ \\
{\rule{0pt}{3.2ex}}{\rule[-1.8ex]{0pt}{0pt}}& $(1)$ & $-\frac{29}{5120}$ \\
{\rule{0pt}{3.2ex}}{\rule[-1.8ex]{0pt}{0pt}}& $(2)$ & $\frac{47}{15360}$ \\
{\rule{0pt}{3.2ex}}{\rule[-1.8ex]{0pt}{0pt}}& $(3)$ & $-\frac{41}{46080}$ \\
{\rule{0pt}{3.2ex}}{\rule[-1.8ex]{0pt}{0pt}}& $(4)$ & $\frac{1}{9216}$ \\
\hline
$(2,0)$  {\rule{0pt}{3.2ex}}{\rule[-1.8ex]{0pt}{0pt}}& $\emptyset$ & $-\mathbf{\frac{1}{384}}$ \\
\hline
\end{tabular}
& 
\begin{tabular}[t]{|c|c|c|}
\hline
$(g,n)$ {\rule{0pt}{3.2ex}}{\rule[-1.8ex]{0pt}{0pt}}& $(k_1,\ldots,k_n)$ & $F_{g,n}[\mathbf{k}]$ \\
\hline
\hline
\multirow{2}{*}{$(1,1)$}  {\rule{0pt}{3.2ex}}{\rule[-1.8ex]{0pt}{0pt}}& $(0)$ & $\mathbf{-\frac{1}{24}}$ \\
 {\rule{0pt}{3.2ex}}{\rule[-1.8ex]{0pt}{0pt}}& $(1)$ & $\frac{1}{48}$ \\
 \hline
\multirow{4}{*}{$(1,2)$}  {\rule{0pt}{3.2ex}}{\rule[-1.8ex]{0pt}{0pt}}& $(0,0)$ & $\mathbf{\frac{1}{32}}$ \\
 {\rule{0pt}{3.2ex}}{\rule[-1.8ex]{0pt}{0pt}}& $(1,0)$& $-\frac{1}{32}$ \\
 {\rule{0pt}{3.2ex}}{\rule[-1.8ex]{0pt}{0pt}}& $(2,0)$ & $-\frac{1}{96}$ \\
 {\rule{0pt}{3.2ex}}{\rule[-1.8ex]{0pt}{0pt}}& $(1,1)$ & $\frac{1}{96}$ \\
\hline
\multirow{7}{*}{$(1,3)$}  {\rule{0pt}{3.2ex}}{\rule[-1.8ex]{0pt}{0pt}}& $(0,0,0)$ & $\mathbf{-\frac{11}{192}}$ \\
 {\rule{0pt}{3.2ex}}{\rule[-1.8ex]{0pt}{0pt}}& $(1,0,0)$ & $\frac{11}{192}$\\
  {\rule{0pt}{3.2ex}}{\rule[-1.8ex]{0pt}{0pt}}& $(2,0,0)$ & $-\frac{5}{192}$ \\
   {\rule{0pt}{3.2ex}}{\rule[-1.8ex]{0pt}{0pt}}& $(1,1,0)$ & $-\frac{1}{24}$ \\
    {\rule{0pt}{3.2ex}}{\rule[-1.8ex]{0pt}{0pt}}& $(3,0,0)$ & $\frac{1}{192}$ \\
     {\rule{0pt}{3.2ex}}{\rule[-1.8ex]{0pt}{0pt}}& $(2,1,0)$ & $\frac{1}{96}$ \\
      {\rule{0pt}{3.2ex}}{\rule[-1.8ex]{0pt}{0pt}}& $(1,1,1)$ & $\frac{1}{96}$  \\
\hline
\end{tabular}
\end{tabular}
	\vspace{2ex}
	\caption{
		For low values of $(g,n)$ we indicate the non-zero values of $F_{g,n}[k_1,\ldots,k_n]$ for $k_1 \geq \cdots \geq k_n \geq 0$ (the others are obtained by symmetry). These coefficients were denoted $F_{g,n}^{{\rm II}}[k_1,\ldots,k_n]$ in Proposition~\ref{secontr}. Inserting the values in bold in \eqref{secondunf} recovers the values of the Masur--Veech volumes given in \cite[Table 11]{TRMVV}. The other values match with the intersection numbers computed from \eqref{fgnseg} via \texttt{admcycles} (\cite{DSvZ}).
	}
	\label{tablak}
\end{table}

\clearpage